\documentclass[a4paper]{article}
\usepackage{amssymb}
\usepackage{graphicx}
\usepackage{hyperref}
\usepackage{latexsym}
\newtheorem{theorem}{Theorem}
\newtheorem{corollary}{Corollary}
\newtheorem{lemma}{Lemma}
\newtheorem{definition}{Definition}

\begin{document}
\title{Some consequences of interpreting the associated logic of the first-order Peano Arithmetic PA finitarily}
\author{\normalsize{Bhupinder Singh Anand}}
\date{\small \textit {Draft of \today. An earlier version of this manuscript is \href{http://arxiv.org/abs/1108.4598}{arXived here}.\footnote{Subject class: LO; MSC: 03B10}}}
\maketitle

\begin{abstract}We show that the classical interpretations of Tarski's inductive definitions actually allow us to define the satisfaction and truth of the quantified formulas of the first-order Peano Arithmetic PA over the domain $N$ of the natural numbers in \textit{two} essentially different ways: (a) in terms of algorithmic verifiabilty; and (b) in terms of algorithmic computability. We show that the classical Standard interpretation $\mathcal{I}_{PA(N,\ Standard)}$ of PA essentially defines the satisfaction and truth of the formulas of the first-order Peano Arithmetic PA in terms of algorithmic verifiability. It is accepted that this classical interpretation---in terms of algorithmic verifiabilty---cannot lay claim to be \textit{finitary}; it does not lead to a finitary justification of the Axiom Schema of Finite Induction of PA from which we may conclude---in an intuitionistically unobjectionable manner---that PA is consistent. We now show that the PA-axioms---including the Axiom Schema of Finite Induction---\textit{are}, however, algorithmically computable finitarily as satisfied / true under the Standard interpretation $\mathcal{I}_{PA(N,\ Standard)}$ of PA; and that the PA rules of inference do preserve algorithmically computable satisfiability / truth finitarily under the Standard interpretation $\mathcal{I}_{PA(N,\ Standard)}$. We conclude that the algorithmically computable PA-formulas \textit{can} provide a finitary interpretation $\mathcal{I}_{PA(N,\ Algorithmic)}$ of PA from which we may classically conclude that PA is consistent in an intuitionistically unobjectionable manner. We define this interpretation, and show that if the associated logic is interpreted finitarily then (i) PA is categorical and (ii) G\"{o}del's Theorem VI holds vacuously in PA since PA is consistent but not $\omega$-consistent. This reflects the fact that PA is $\omega$-consistent if, and only if, Aristotle's particularisation is presumed to always hold under any interpretation of the associated logic; and that the standard interpretation  of PA is a model of PA if, and only if, PA is $\omega$-consistent. 

\vspace{+1ex}
\noindent \scriptsize \textbf{Keywords} Algorithmic computability, algorithmic verifiability, Aristotle's particularisation, consistency, first-order, $\omega$-consistency, Peano Arithmetic PA, satisfaction, soundness, standard interpretation, Tarski.
\end{abstract}

\section{Introduction}
\label{sec:1}

In a recent paper\footnote{\cite{An12}.}, `Evidence-Based Interpretations of PA', presented at the Symposium on Computational Philosophy at the AISB/IACAP 2012 World Congress, Birmingham, we showed first that---in addition to the classically defined Standard interpretation $\mathcal{I}_{PA(N,\ Standard)}$ of the first order Peano Arithmetic PA over the domain $N$ of the natural numbers---Tarski's classical definitions of the satisfaction and truth of the formulas of a formal language under an interpretation admit two evidence-based interpretations of PA under the standard first order logic FOL:

\begin{itemize}
\item
An Instantiational interpretation $\mathcal{I}_{PA(\mathbb{N},\ Instantiational)}$ of PA over the domain $\mathbb{N}$ of the PA numerals; and

\item
An Algorithmic interpretation $\mathcal{I}_{PA(N,\ Algorithmic)}$ of PA over the domain $N$ of the natural numbers.
\end{itemize}

\noindent We then showed that the Instantiational interpretation $\mathcal{I}_{PA(\mathbb {N},\ Instantiational)}$ of PA is sound if, and only if, the Standard interpretation $\mathcal{I}_{PA(N,\ Standard)}$ of PA is sound; where we defined an interpretation of PA as sound if, and only if:

\begin{itemize}
\item
The axioms of PA are true under the interpretation; and

\item
The PA rules of inference preserve such truth.
\end{itemize}

\noindent We further showed that:

\begin{itemize}
\item
The axioms of PA are true under the Algorithmic interpretation $\mathcal{I}_{PA(N,\ Algo-}$ $_{rithmic)}$ of PA; and

\item
The PA rules of inference preserve such truth.
\end{itemize}

\noindent We concluded that:

\begin{itemize}
\item
The Algorithmic interpretation $\mathcal{I}_{PA(N,\ Algorithmic)}$ of PA is sound; and

\item
PA is consistent.
\end{itemize}

\subsection{The philosophical question addressed by our investigation}

Our investigation sought to address the question:

\begin{itemize}
\item
Is there any objective evidence to justify the acceptance of arithmetical propositions as `true' on the grounds that such `truth' is self-evident?
\end{itemize}

\noindent We noted, for instance, that conventional wisdom follows Tarski's inductive definitions of the `satisfiability' and `truth' of the formulas of a formal language such as PA under an interpretation when it implicitly and, as we shall show\footnote{In Theorem \ref{sec:6.3.cor.3}.}, non-finitarily\footnote{We use the term `finitary' essentially in the broader sense detailed by David Hilbert in \cite{Hi25}. See also the Stanford Encyclopedia of Philosophy article on `Hilbert's Program: The Finitary Point of View' at \href{http://plato.stanford.edu/entries/hilbert-program/\#2}{\scriptsize http://plato.stanford.edu/entries/hilbert-program/\#2}.} holds that:

\begin{itemize}
\item
The Standard interpretation $\mathcal{I}_{PA(N,\ Standard)}$ of PA over the domain $N$ of the natural numbers is sound if the standard interpretation of FOL is sound.
\end{itemize}

\noindent In other words, conventional wisdom holds it as self-evident that---even though an infinite process is implicit in their decidability:

\begin{itemize}
\item
The denumerable atomic formulas of PA can be \textit{assumed} as decidable under the Standard interpretation $\mathcal{I}_{PA(N,\ Standard)}$ of PA;
    
\vspace{+1ex}
\item
The denumerable PA axioms can be \textit{assumed} to interpret as true under the Standard interpretation $\mathcal{I}_{PA(N,\ Standard)}$ of PA;
     
\vspace{+1ex}
\item
The PA rules of inference can be \textit{assumed} to preserve truth under the Standard interpretation $\mathcal{I}_{PA(N,\ Standard)}$ of PA.
\end{itemize}

\noindent What this means is that conventional wisdom also holds it as self-evident under the standard interpretation of FOL that---even though an infinite process is implicit in the decidability:

\begin{itemize}
\item
The formula $[(\forall x)F(x)]$ is decidable under the Standard interpretation $\mathcal{I}_{PA(N,}$ $_{Standard)}$ of PA.\footnote{We note that, as emphasised by Edward Nelson in \cite{Ne00}, the assumption of the unqualified decidability of quantified formulas under Tarski's definitions has been a matter of controversy. Theorem \ref{sec:6.3.cor.3} now shows that the assumption is indeed untenable.}

\vspace{+1ex}
\item
If the formula $[(\exists x)F(x)]$ is true under the Standard interpretation $\mathcal{I}_{PA(N,}$ $_{Standard)}$ of PA, then there must exist some numeral $[n]$ for which the formula $[F(n)]$ is true under the interpretation.\footnote{In other words Aristotle's particularisation (Definition \ref{sec:A.def.1}) is valid---without qualification---over $N$.} 
\end{itemize}

\noindent We also noted that---unless we assume that PA is $\omega$-consistent---we cannot conclude by FOL that:

\begin{itemize}
\item
If the formula $[(\exists x)F(x)]$ is provable in PA, then there must exist some numeral $[n]$ for which the formula $[F(n)]$ is provable in PA. 
\end{itemize}

\subsection{What differentiates our approach?}
\label{sec:1.03.1}

Our approach to the above investigation can be differentiated by noting first that, in comparison, conventional wisdom---essentially following David Hilbert\footnote{In a 1925 address (\cite{Hi25}) Hilbert had shown that the axiomatisation $\mathcal{L}_{\varepsilon}\) of classical Aristotlean predicate logic proposed by him as a formal first-order $\varepsilon\)-predicate calculus (detailed in \cite{Hi27}, pp.465-466) in which he used a primitive choice-function (\cite{Hi25}, p.382) symbol, `$\varepsilon\)', for defining the quantifiers `$\forall\)' and `$\exists\)' would adequately express---and yield, under a suitable interpretation---Aristotle's logic of predicates if the $\varepsilon\)-function was interpreted to yield Aristotlean particularisation (\cite{Hi25}, pp.382-383; \cite{Hi27}, p.466(1)).}---can be labelled `theistic' in that it implicitly assumes both that:

\begin{itemize}
\item
The standard first order logic FOL is consistent;
\end{itemize}

\noindent and that:

\begin{itemize}
\item
The standard interpretation of FOL is sound.
\end{itemize}

\begin{quote}
\noindent The significance of the label `theistic' is that conventional wisdom tacitly believes that Aristotle's particularisation\footnote{Definition \ref{sec:A.def.1}.} remains valid---without qualification---even over infinite domains; a belief that is not unequivocally self-evident, but must be appealed to as an article of faith.
\end{quote}

\noindent We note second that, in sharp contrast, constructive approaches to mathematics---such as Intuitionism---can be labelled `atheistic' since they deny both that:

\begin{itemize}
\item
FOL is consistent (since they deny the Law of The Excluded Middle\footnote{``The formula $\forall x(A(x)\vee \neg A(x))$ is classically provable, and hence under classical interpretation true. But it is unrealizable. So if realizability is accepted as a necessary condition for intuitionistic truth, it is untrue intuitionistically, and therefore unprovable not only in the present intuitionistic formal system, but by any intuitionistic methods whatsoever". \cite{Kl52}, p.513.}.);
\end{itemize}

\noindent and that:

\begin{itemize}
\item
The standard interpretation of FOL is sound (since they deny Aristotle's particularisation).
\end{itemize}

\begin{quote}
\noindent The significance of the label `atheistic' is that whereas constructive approaches to mathematics deny the faith-based belief in the validity of Aristotle's particularisation---without qualification---over infinite domains, their denial of the Law of the Excluded Middle is itself a belief---in the inconsistency of FOL---that is also not unequivocally self-evident, and must also be appealed to as an article of faith\footnote{Although Brouwer's explicitly stated objection appeared to be to the Law of the Excluded Middle as expressed and interpreted at the time (\cite{Br23}, p.335-336; \cite{Kl52}, p.47; \cite{Hi27}, p.475), some of Kleene's remarks (\cite{Kl52}, p.49), some of Hilbert's remarks (for instance in \cite{Hi27}, p.474) and, more particularly, Kolmogorov's remarks (in \cite{Ko25}, fn. p.419; p.432) suggest that the intent of Brouwer's fundamental objection can also be viewed today as being limited only to the yet prevailing belief---as an article of faith---that the validity of Aristotle's particularisation can be extended without qualification to infinite domains.}.
\end{quote}

\noindent In our investigation, however, we follow what may be labelled an `agnostic' approach by noting that although, if Aristotle's particularisation holds in an interpretation then the Law of the Excluded Middle must also hold in the interpretation, the converse is not true.

\vspace{+1ex}
\noindent We thus follow a middle path by explicitly assuming that:

\begin{itemize}
\item
FOL is consistent;
\end{itemize}

\noindent and explicitly state when an argument appeals to the postulation that:

\begin{itemize}
\item
The standard interpretation of FOL is sound.
\end{itemize}

\begin{quote}
\noindent The significance of the label `agnostic' is that we neither hold FOL to be inconsistent, nor hold that Aristotle's particularisation can be applied---without qualification---over infinite domains.
\end{quote}

\section{Overview}
\label{sec:2.1}

In this paper we revisit the arguments of \cite{An12} and consider some consequences.

\vspace{+1ex}
\noindent Specifically, we first define what it means for a formula of an arithmetical language such as the first order Peano Arithmetic PA to be:

\begin{quote}
(i) Algorithmically verifiable (Definition \ref{sec:1.02.def.1});

\vspace{+.5ex}
(ii) Algorithmically computable (Definition \ref{sec:1.02.def.2}).
\end{quote}

\noindent under an interpretation.

\vspace{+1ex}
\noindent We then show that:

\begin{quote}
(a) The PA-formulas are decidable under the standard interpretation of PA if, and only if, they are algorithmically verifiable under the interpretation (Corollary \ref{sec:5.5.cor.2});

\begin{quote}
\footnotesize Although the standard interpretation is believed to define a model of PA, the definition cannot claim to be finitary since it does not lead to a finitary justification of the Axiom Schema of (finite) Induction of PA from which we may conclude---in an intuitionistically unobjectionable manner---that PA is consistent\footnote{The possibility/impossibility of such justification was the subject of the famous Poincar\'{e}-Hilbert debate. See \cite{Hi27}, p.472; also \cite{Br13}, p.59; \cite{We27}, p.482; \cite{Pa71}, p.502-503.}. We note further that Gerhard Gentzen's `constructive'\footnote{In the sense highlighted by Elliott Mendelson in \cite{Me64}, p.261.} consistency proof for formal number theory\footnote{cf. \cite{Me64}, p258.} is debatably finitary\footnote{See for instance \href{http://en.wikipedia.org/wiki/Hilbert's_program}{http://en.wikipedia.org/wiki/Hilbert's\_program}.}, since it involves a Rule of Infinite Induction that admits appeal to the well-ordering property of transfinite ordinals; and we show in Section \ref{sec:10}, Appendix E that we cannot introduce a transfinite ordinal into any model of PA without inviting inconsistency.
\end{quote}

\vspace{+,5ex}
(b) The PA-axioms are algorithmically computable as satisfied / true under the standard interpretation of PA (Lemmas \ref{sec:6.2.lem.1} and \ref{sec:6.2.lem.2});

\vspace{+.5ex}
(c) Generalisation and Modus Ponens preserve algorithmically computable truth under the standard interpretation of PA (Lemmas \ref{sec:6.2.lem.3} and \ref{sec:6.2.lem.4});

\vspace{+.5ex}
(d) The provable PA-formulas are precisely the ones that are algorithmically computable as satisfied / true under the standard interpretation of PA (Theorem \ref{sec:6.2.lem.5}).
\end{quote}

\noindent We conclude that the algorithmically computable PA-formulas can provide a sound---in the sense of Definition \ref{sec:A.def.10}---finitary interpretation of PA (Theorem \ref{sec:6.2.thm.1}).

\vspace{+1ex}
\noindent We note that PA is $\omega$-consistent if, and only if, Aristotle's particularisation (Definition \ref{sec:A.def.1}) is presumed to always hold under any interpretation of the associated logic (Section \ref{sec:4}, Appendix A). 

\vspace{+1ex}
\noindent We then show that if classical first-order logic is interpreted finitarily (Section \ref{sec:8}, Appendix C) without the presumption that Aristotle's particularisation necessarily holds under the interpretation, then we may conclude that:

\begin{quote}

(f) PA is consistent (Theorem \ref{sec:6.2.thm.2});

\vspace{+1ex}
(g) PA is categorical (Corollary \ref{sec:6.3.thm.1.cor.1});

\vspace{+1ex}
(h) PA is not $\omega$-consistent (Corollary \ref{sec:6.3.cor.2});

\vspace{+1ex}
(i) the standard interpretation of PA is not sound\footnote{In the sense of Definition \ref{sec:A.def.10}.}, and does not yield a model for PA (Corollary \ref{sec:6.3.cor.3}).
\end{quote}

\subsection{Notation, Definitions and Comments}
\label{sec:A}

\textbf{Comments} We have taken some liberty in emphasising standard definitions selectively, and interspersing our arguments liberally with comments and references, generally of a foundational nature. These are intended to reflect our underlying thesis that essentially arithmetical problems appear more natural when expressed---and viewed---within the perspective of an interpretation of PA that appeals to the \textit{evidence} provided by a deterministic algorithm along the lines suggested in Section \ref{sec:5.4.0}; a perspective that, by its very nature, cannot appeal implicitly to transfinite concepts.

\begin{quote}
\label{sec:A.qtn}
\footnotesize \textbf{Evidence} ``It is by now folklore \ldots that one can view the \textit{values} of a simple functional language as specifying \textit{evidence} for propositions in a constructive logic \ldots"\footnote{\cite{Mu91}.}.
\end{quote}

\vspace{+1ex}
\noindent \textbf{Notation} We use square brackets to indicate that the contents represent a symbol or a formula---of a formal theory---generally assumed to be well-formed unless otherwise indicated by the context.\footnote{In other words, expressions inside the square brackets are to be only viewed syntactically as juxtaposition of symbols that are to be formed and manipulated upon strictly in accordance with specific rules for such formation and manipulation---in the manner of a mechanical or electronic device---without any regards to what the symbolism might represent semantically under an interpretation that gives them meaning.} We use an asterisk to indicate that the associated expression is to be interpreted semantically with respect to some well-defined interpretation.

\begin{definition}
\label{sec:A.def.1}
\noindent \textbf{Aristotle's particularisation} This holds that from a meta-assertion such as:

\begin{quote}
`It is not the case that: For any given $x$, $P^{*}(x)$ does not hold',
\end{quote}

\noindent usually denoted symbolically by `$\neg(\forall x)\neg P^{*}(x)$', we may always validly infer in the classical, Aristotlean, logic of predicates\footnote{\cite{HA28}, pp.58-59.} that:

\begin{quote}
`There exists an unspecified $x$ such that $P^{*}(x)$ holds',
\end{quote}

\noindent usually denoted symbolically by `$(\exists x)P^{*}(x)$'.
\end{definition}

\begin{quote}
\footnotesize
\textbf{The significance of Aristotle's particularisation for the first-order predicate calculus:} We note that in a formal language the formula `$[(\exists x)P(x)]$' is an abbreviation for the formula `$[\neg(\forall x)\neg P(x)]$'. The commonly accepted interpretation of this formula---and a fundamental tenet of classical logic unrestrictedly adopted as intuitively obvious by standard literature\footnote{See \cite{Hi25}, p.382; \cite{HA28}, p.48; \cite{Sk28}, p.515; \cite{Go31}, p.32.; \cite{Kl52}, p.169; \cite{Ro53}, p.90; \cite{BF58}, p.46; \cite{Be59}, pp.178 \& 218; \cite{Su60}, p.3; \cite{Wa63}, p.314-315; \cite{Qu63}, pp.12-13; \cite{Kn63}, p.60; \cite{Co66}, p.4; \cite{Me64}, p.52(ii); \cite{Nv64}, p.92; \cite{Li64}, p.33; \cite{Sh67}, p.13; \cite{Da82}, p.xxv; \cite{Rg87}, p.xvii; \cite{EC89}, p.174; \cite{Mu91}; \cite{Sm92}, p.18, Ex.3; \cite{BBJ03}, p.102.} that seeks to build upon the formal first-order predicate calculus---tacitly appeals to Aristotlean particularisation.

However, L. E. J. Brouwer had noted in his seminal 1908 paper on the unreliability of logical principles\footnote{\cite{Br08}.} that the commonly accepted interpretation of this formula is ambiguous if interpretation is intended over an infinite domain.

Brouwer essentially argued that, even supposing the formula `$[P(x)]$' of a formal Arithmetical language interprets as an arithmetical relation denoted by `$P^{*}(x)$', and the formula `$[\neg(\forall x)\neg P(x)]$' as the arithmetical proposition denoted by `$\neg(\forall x)\neg P^{*}(x)$', the formula `$[(\exists x)P(x)]$' need not interpret as the arithmetical proposition denoted by the usual abbreviation `$(\exists x)P^{*}(x)$'; and that such postulation is invalid as a general logical principle in the absence of a means for constructing some putative object $a$ for which the proposition $P^{*}(a)$ holds in the domain of the interpretation.

Hence we shall follow the convention that the assumption that `$(\exists x)P^{*}(x)$' is the intended interpretation of the formula `$[(\exists x)P(x)]$'---which is essentially the assumption that Aristotle's particularisation holds over the domain of the interpretation---must always be explicit.

\vspace{+1ex}
\noindent \textbf{The significance of Aristotle's particularisation for PA:} In order to avoid intuitionistic objections to his reasoning, Kurt G\"{o}del introduced the syntactic property of $\omega$-consistency\footnote{The significance of $\omega$-consistency for the formal system PA is highlighted in Section \ref{sec:4}, Appendix A.} as an explicit assumption in his formal reasoning in his seminal 1931 paper on formally undecidable arithmetical propositions\footnote{\cite{Go31}, p.23 and p.28.}.

G\"{o}del explained at some length\footnote{In his introduction on p.9 of \cite{Go31}.} that his reasons for introducing $\omega$-consistency explicitly was to avoid appealing to the semantic concept of classical arithmetical truth in Aristotle's logic of predicates (which presumes Aristotle's particularisation).

We show in Section \ref{sec:4}, Appendix A that the two concepts are meta-mathematically equivalent in the sense that, if PA is consistent, then PA is $\omega$-consistent if, and only if, Aristotle's particularisation holds under the standard interpretation of PA.
\end{quote}

\begin{definition}
\label{sec:A.def.2}
\noindent \textbf{The structure of the natural numbers:} \{$N$ (\textit{the set of natural numbers}); $=$ (\textit{equality}); $'$ (\textit{the successor function}); $+$ (\textit{the addition function}); $ \ast $ (\textit{the product function}); $0$ (\textit{the null} \textit{element})\}.
\end{definition}

\begin{definition}
\label{sec:A.def.3}
\noindent \textbf{The axioms of first-order Peano Arithmetic (PA)}
\end{definition}

\begin{tabbing}
\textbf{PA$_{1}$} \= $[(x_{1} = x_{2}) \rightarrow ((x_{1} = x_{3}) \rightarrow (x_{2} = x_{3}))]$; \\

\textbf{PA$_{2}$} \> $[(x_{1} = x_{2}) \rightarrow (x_{1}^{\prime} = x_{2}^{\prime})]$; \\

\textbf{PA$_{3}$} \> $[0 \neq x_{1}^{\prime}]$; \\

\textbf{PA$_{4}$} \> $[(x_{1}^{\prime} = x_{2}^{\prime}) \rightarrow (x_{1} = x_{2})]$; \\

\textbf{PA$_{5}$} \> $[( x_{1} + 0) = x_{1}]$; \\

\textbf{PA$_{6}$} \> $[(x_{1} + x_{2}^{\prime}) = (x_{1} + x_{2})^{\prime}]$; \\

\textbf{PA$_{7}$} \> $[( x_{1} \star 0) = 0]$; \\

\textbf{PA$_{8}$} \> $[( x_{1} \star x_{2}^{\prime}) = ((x_{1} \star x_{2}) + x_{1})]$; \\

\textbf{PA$_{9}$} \> For any well-formed formula $[F(x)]$ of PA: \\

\> $[F(0) \rightarrow (((\forall x)(F(x) \rightarrow F(x^{\prime}))) \rightarrow (\forall x)F(x))]$.
\end{tabbing}

\begin{definition}
\label{sec:A.def.4}
\noindent \textbf{Generalisation in PA} If $[A]$ is PA-provable, then so is $[(\forall x)A]$.
\end{definition}

\begin{definition}
\label{sec:A.def.5}
\noindent \textbf{Modus Ponens in PA} If $[A]$ and $[A \rightarrow B]$ are PA-provable, then so is $[B]$.
\end{definition}

\begin{definition}
\label{sec:A.def.6}
\noindent \textbf{Standard interpretation of PA} The standard interpretation $\mathcal{I}_{PA(N,\ Standard)}$ of PA over the structure $N$ is the one in which the logical constants have their `usual' interpretations\footnote{These are expressed formally in Section \ref{sec:7}, Appendix B, and essentially follow the definitions in \cite{Me64}, p.49.} in Aristotle's logic of predicates (which subsumes Aristotle's particularisation), and\footnote{See \cite{Me64}, p.107.}:

\begin{tabbing}
(a) \= the set of non-negative integers is the domain; \\

(b) \> the symbol [0] interprets as the integer 0; \\

(c) \> the symbol $[']$ interprets as the successor operation (addition of 1); \\

(d) \> the symbols $[+]$ and $[*]$ interpret as ordinary addition and multiplication; \\

(e) \> the symbol $[=]$ interprets as the identity relation.
\end{tabbing}
\end{definition}

\begin{definition}
\label{sec:A.def.7}
\noindent \textbf{Simple consistency:} A formal system S is simply consistent if, and only if, there is no S-formula $[F(x)]$ for which both $[(\forall x)F(x)]$ and $[\neg(\forall x)F(x)]$ are S-provable.
\end{definition}

\begin{definition}
\label{sec:A.def.8}
\noindent \textbf{$\omega$-consistency:} A formal system S is $\omega$-consistent if, and only if, there is no S-formula $[F(x)]$ for which, first, $[\neg(\forall x)F(x)]$ is S-provable and, second, $[F(a)]$ is S-provable for any given S-term $[a]$.
\end{definition}

\begin{definition}
\label{sec:A.def.9}
\noindent \textbf{Soundness (formal system - non-standard):}  A formal system S is sound under an interpretation $\mathcal{I}_{S}$ with respect to a domain $\mathbb{D}$ if, and only if, every theorem $[T]$ of S translates as `$[T]$ is true under $\mathcal{I}_{S}$ in $\mathbb{D}$'.
\end{definition}

\begin{definition}
\label{sec:A.def.10}
\noindent \textbf{Soundness (interpretation - non-standard):}  An interpretation $\mathcal{I}_{S}$ of a formal system S is sound with respect to a domain $\mathbb{D}$ if, and only if, S is sound under the interpretation $\mathcal{\mathcal{I}_{S}}$ over the domain $\mathbb{D}$.
\end{definition}

\begin{quote}
\scriptsize{\textbf{Soundness in classical logic:} In classical logic, a formal system $S$ is sometimes defined as `sound' if, and only if, it has an interpretation; and an interpretation is defined as the assignment of meanings to the symbols, and truth-values to the sentences, of the formal system. Moreover, any such interpretation is defined as a model\footnote{We follow the definition in \cite{Me64}, p.51.} of the formal system. This definition suffers, however, from an implicit circularity: the formal logic $L$ underlying any interpretation of $S$ is implicitly assumed to be `sound'. The above definitions seek to avoid this implicit circularity by delinking the defined `soundness' of a formal system under an interpretation from the implicit `soundness' of the formal logic underlying the interpretation. This admits the case where, even if $L_{1}$ and $L_{2}$ are implicitly assumed to be sound, $S+L_{1}$ is sound, but $S+L_{2}$ is not. Moreover, an interpretation of $S$ is now a model for $S$ if, and only if, it is sound.\footnote{My thanks to Professor Rohit Parikh for highlighting the need for making such a distinction explicit.}}
\end{quote}

\begin{definition}
\label{sec:A.def.11}
\noindent \textbf{Categoricity:} A formal system S is categorical if, and only if, it has a sound\footnote{In the sense of Definitions \ref{sec:A.def.9} and \ref{sec:A.def.10}.} interpretation and any two sound interpretations of S are isomorphic.\footnote{Compare \cite{Me64}, p.91.}
\end{definition}

\section{Interpretation of an arithmetical language in terms of the computations of a simple functional language}
\label{sec:1.03}

We begin by noting that we can, in principle, define\footnote{Formal definitions are given in Section \ref{sec:5.4.0}.} the classical `satisfaction' and `truth' of the formulas of a first order arithmetical language, such as PA, \textit{verifiably} under an interpretation using as \textit{evidence}\footnote{\cite{Mu91}.} the computations of a simple functional language.

\vspace{+1ex}
\noindent Such definitions follow straightforwardly for the atomic formulas of the language (i.e., those without the logical constants that correspond to `negation', `conjunction', `implication' and `quantification') from the standard definition of a simple functional language\footnote{Such as, for instance, that of a deterministic Turing machine (\cite{Me64}, pp.229-231) based essentially on Alan Turing's seminal 1936 paper on computable numbers (\cite{Tu36}).}.

\vspace{+1ex}
\noindent Moreover, following Alfred Tarski's seminal 1933 paper on the the concept of truth in the languages of the deductive sciences\footnote{\cite{Ta33}.}, the classical `satisfaction' and `truth' of those formulas of a first-order language which contain logical constants can be inductively defined, under an interpretation, in terms of the `satisfaction' and `truth' of the interpretations of only the atomic formulas of the language.

\vspace{+1ex}
\noindent Hence, classically, the `satisfaction' and `truth' of those formulas of an arithmetical language such as PA which contain logical constants can, in principle, also be defined verifiably under an interpretation using as evidence the computations of a simple functional language.

\vspace{+1ex}
\noindent We show in Section \ref{sec:5.4.0} that this is indeed the case for PA under the standard interpretation $\mathcal{I}_{PA(N,\ Standard)}$, when this is explicitly defined as in Section \ref{sec:5.4.0.0}.

\vspace{+1ex}
\noindent We show, moreover, that we can further define `algorithmic truth' and `algorithmic falsehood' finitarily under $\mathcal{I}_{PA(N,\ Standard)}$ such that the PA axioms interpret as always algorithmically true, and the rules of inference preserve algorithmic truth, over the domain $N$ of the natural numbers.

\begin{quote}
\footnotesize
\textbf{Significance of `algorithmic truth':} The \textit{algorithmically} true propositions of $N$ under $\mathcal{I}_{PA(N,\ Standard)}$ are, moreover, a proper subset\footnote{This follows immediately from Corollary \ref{sec:6.3.cor.1}.} of the \textit{verifiably} true propositions of $N$ under $\mathcal{I}_{PA(N,\ Standard)}$; they suggest a possible finitary model of PA that establishes the consistency of PA constructively.
\end{quote}

\subsection{The definitions of `algorithmic truth' and `algorithmic falsehood' under $\mathcal{I}_{PA(N,\ Standard)}$ are not symmetric with respect to `truth' and `falsehood' under $\mathcal{I}_{PA(N,\ Standard)}$}
\label{sec:1.03.1}

However, the definitions of `algorithmic truth' and `algorithmic falsehood' under $\mathcal{I}_{PA(N,\ Standard)}$ are not symmetric with respect to classical (verifiable) `truth' and `falsehood' under $\mathcal{I}_{PA(N,\ Standard)}$.

\vspace{+1ex}
\noindent For instance, if a formula $[(\forall x)F(x)]$ of an arithmetic is algorithmically true under an interpretation (such as $\mathcal{I}_{PA(N,\ Standard)}$), then we may conclude that there is a deterministic algorithm that, for any given numeral $[a]$, provides evidence that the formula $[F(a)]$ is algorithmically true under the interpretation.

\vspace{+1ex}
\noindent In other words, there is a deterministic algorithm that provides evidence that the interpretation $F^{*}(a)$ of $[F(a)]$ \textit{holds} in $N$ for any given natural number $a$.

\begin{quote}
\footnotesize
\textbf{Defining the term `hold':} We define the term `hold'---when used in connection with an interpretation of a formal language and, more specifically, with reference to the computations of a simple functional language associated with the atomic formulas of the language---explicitly in Section \ref{sec:5.4.0}; the aim being to avoid appealing to the classically subjective (and existential) connotation implicitly associated with the term under an implicitly defined standard interpretation of an arithmetic\footnote{As, for instance, in \cite{Go31}.}.
\end{quote}

\noindent However, if a formula $[(\forall x)F(x)]$ of an arithmetic is algorithmically false under an interpretation, then we can only conclude that there is no deterministic algorithm that, for any given natural number $a$, can provide evidence whether the interpretation $F^{*}(a)$ holds or not in $N$ .

\vspace{+1ex}
\noindent We cannot therefore conclude that there is a numeral $[a]$ such that the formula $[F(a)]$ is algorithmically false under the interpretation; nor can we conclude that there is a natural number $b$ such that $F^{*}(b)$ does not hold in $N$.

\vspace{+1ex}
\noindent Such a conclusion would require:

\begin{quote}
(i) either some additional evidence that will verify for some assignment of numerical values to the free variables of $[F]$ that the corresponding interpretation $F^{*}$ does not hold\footnote{Essentially reflecting Brouwer's objection to the assumption of Aristotle's particularisation over an infinite domain.};

\vspace{+1ex}
(ii) or the additional assumption that either Aristotle's particularisation\footnote{Definition \ref{sec:A.def.1}.} holds over the domain of the interpretation (as is implicitly presumed under the standard interpretation of PA) or that the arithmetic is $\omega$-consistent\footnote{An assumption explicitly introduced by G\"{o}del in \cite{Go31}.}.
\end{quote}

\section{Defining algorithmic verifiability and algorithmic computability}
\label{sec:1.02}

The asymmetry of Section \ref{sec:1.03.1} suggests\footnote{My thanks to Dr. Chaitanya H. Mehta for advising that the focus of this investigation should be the distinction between these two concepts.} the following two concepts\footnote{When dealing with infinite processes, the distinction sought to be made between algorithmically verifiable formulas and algorithmically computable formulas can be viewed as reflecting in number theory the analogous distinction that is made in analysis between, for instance, continuous functions (\cite{Ru53}, p.65, \S 4.5) and uniformly continuous functions (\cite{Ru53}, p.65, \S 4.13); or that between convergent sequences (\cite{Ru53}, p.65, \S 7.1) and uniformly convergent sequences (\cite{Ru53}, p.65, \S 7.7).}: 

\begin{definition}
\label{sec:1.02.def.1}
\textbf{Algorithmic verifiability: single variable} 

\vspace{+1ex}
\noindent An arithmetical relational formula $[F(x)]$ is algorithmically verifiable under an interpretation if, and only if, for any given numeral $[n]$, we can define a deterministic algorithm $AL_{n}$ which provides objective evidence for deciding the truth/falsity of each proposition in the finite sequence $\{[F(1)], [F(2)], \ldots, [F(n)]\}$ under the interpretation.
\end{definition}

\noindent \textbf{\textit{\scriptsize Example:}} {\scriptsize Since any real number is definable as the limit of a Cauchy sequence of rational numbers:}

\begin{itemize}
\scriptsize
\item
Let $[R(n)]$ denote the $n^{th}$ digit in the decimal expression of the real number $R$ in binary notation.

\vspace{+1ex}
\item
Then, for any given natural number $n$, there is a deterministic algorithm $AL_{n}$ that will decide the truth/falsity of each proposition in the sequence $\{[R(1)=0], [R(2)=0], \ldots, [R(n)=0]\}$.

\vspace{+1ex}
\item
Hence $[R(x)=0]$ is algorithmically verifiable.
\end{itemize}

\begin{definition}
\label{sec:1.02.def.2}
\textbf{Algorithmic computability: single variable} 

\vspace{+1ex}
\noindent An arithmetical relational formula $[F(x)]$ is algorithmically computable under an interpretation if, and only if, we can define a deterministic algorithm $AL$ that provides objective evidence for deciding the truth/falsity of each proposition in the denumerable sequence $\{[F(1), [F(2)], \ldots\}]$ under the interpretation.
\end{definition}

\noindent We note that although every algorithmically computable formula with a single variable is algorithmically verifiable, the converse is not true.

\vspace{+1ex}
\noindent \textbf{\textit{\scriptsize Example:}} {\scriptsize Since it follows from Alan Turing's Halting argument\footnote{\cite{Tu36}, p.132, \S 8.} that there are algorithmically uncomputable real numbers:}

\begin{itemize}
\scriptsize
\item
Let $[R(n)]$ denote the $n^{th}$ digit in the decimal expression of an algorithmically uncomputable real number $R$ in binary notation.

\vspace{+1ex}
\item
Then, for any given natural number $n$, there is a deterministic algorithm $AL_{n}$ that will decide the truth/falsity of each proposition in the sequence $\{[R(1)=0], [R(2)=0], \ldots, [R(n)=0]\}$.

\vspace{+1ex}
\item
However, there is no deterministic algorithm $AL$ that will decide the truth/falsity of each proposition in the denumerable sequence $\{[R(1)=0], [R(2)=0], \ldots\}$.

\vspace{+1ex}
\item
Hence the relational formula $[R(x)=0]$ is algorithmically verifiable but not algorithmically computable.
\end{itemize}

\noindent We note that we can generalise Definition \ref{sec:1.02.def.1} to:

\begin{definition}
\label{sec:1.02.def.1.1}
\textbf{Algorithmic verifiability} 

\vspace{+1ex}
\noindent An arithmetical relational formula $[F(x_{1}, x_{2}, \ldots, x_{k})]$ is algorithmically verifiable under an interpretation if, and only if, for any given sequence of numerals $[a_{1}, a_{2}, \ldots, a_{k}]$, we can define a deterministic algorithm $AL_{a}$ which provides objective evidence for deciding the truth/falsity of the formula $[F(a_{1}, a_{2}, \ldots, a_{k})]$ under the interpretation.
\end{definition}

\noindent We show in Section \ref{sec:5.4.0} that the `algorithmic verifiability' of the formulas of a formal language which contain logical constants can be inductively defined under an interpretation in terms of the `algorithmic verifiability' of the interpretations of the atomic formulas of the language; further, that the PA-formulas are decidable under the standard interpretation of PA if, and only if, they are algorithmically verifiable under the interpretation (Corollary \ref{sec:5.5.cor.2}).

\vspace{+1ex}
\noindent We can similarly generalise Definition \ref{sec:1.02.def.2} to:

\begin{definition}
\label{sec:1.02.def.2.1}
\textbf{Algorithmic computability} 

\vspace{+1ex}
\noindent An arithmetical relational formula $[F(x_{1}, x_{2}, \ldots, x_{k})]$ is algorithmically computable under an interpretation if, and only if, we can define a deterministic algorithm $AL$ that, for any given sequence of numerals $[a_{1}, a_{2}, \ldots, a_{k}]$, provides objective evidence for deciding the truth/falsity of the proposition $[F(a_{1}, a_{2}, \ldots, a_{k})]$ under the interpretation.
\end{definition}

\noindent We show in Section \ref{sec:5.4.0} that the `algorithmic computability' of the formulas of a formal language which contain logical constants can also be inductively defined under an interpretation in terms of the `algorithmic computability' of the interpretations of the atomic formulas of the language; further, that the PA-formulas are decidable under an algorithmic interpretation of PA if, and only if, they are algorithmically computable under the interpretation \label{sec:6.3.thm.1}.

\vspace{+1ex}
\noindent We now show that the above concepts are well-defined under the standard interpretation of PA.

\section{The \textit{implicit} Satisfaction condition in Tarski's inductive assignment of truth-values under an interpretation}
\label{sec:5.4.0}

We first consider the significance of the \textit{implicit} Satisfaction condition in Tarski's inductive assignment of truth-values under an interpretation.

\vspace{+1ex}
\noindent We note that---essentially following standard expositions\footnote{cf. \cite{Me64}, p.51.} of Tarski's inductive definitions on the `satisfiability' and `truth' of the formulas of a formal language under an interpretation---we can define:

\begin{definition}
\label{sec:5.def.1}
If $[A]$ is an atomic formula $[A(x_{1}, x_{2}, \ldots, x_{n})]$ of a formal language S, then the denumerable sequence $(a_{1}, a_{2}, \ldots)$ in the domain $\mathbb{D}$ of an interpretation $\mathcal{I}_{S(\mathbb{D})}$ of S satisfies $[A]$ if, and only if:

\begin{quote}
(i) $[A(x_{1}, x_{2}, \ldots, x_{n})]$ interprets under $\mathcal{I}_{S(\mathbb{D})}$ as a unique relation $A^{*}(x_{1}, x_{2},$ $\ldots, x_{n})$ in $\mathbb{D}$ for any witness $\mathcal{W}_{\mathbb{D}}$ of $\mathbb{D}$;

\vspace{+1ex}
(ii) there is a Satisfaction Method, SM($\mathcal{I}_{S(\mathbb{D})}$) that provides objective \textit{evidence}\footnote{In the sense of \cite{Mu91}.} by which any witness $\mathcal{W}_{\mathbb{D}}$ of $\mathbb{D}$ can objectively \textbf{define} for any atomic formula $[A(x_{1}, x_{2}, \ldots, x_{n})]$ of S, and any given denumerable sequence $(b_{1}, b_{2}, \ldots)$ of $\mathbb{D}$, whether the proposition $A^{*}(b_{1}, b_{2}, \ldots, b_{n})$ holds or not in $\mathbb{D}$;

\vspace{+1ex}
(iii) $A^{*}(a_{1}, a_{2}, \ldots, a_{n})$ holds in $\mathbb{D}$ for any $\mathcal{W}_{\mathbb{D}}$.
\end{quote}
\end{definition}

\begin {quote}
\footnotesize
\textbf{Witness:} From a finitary perspective, the existence of a `witness' as in (i) above is implicit in the usual expositions of Tarski's definitions.

\vspace{+1ex}
\textbf{Satisfaction Method:} From a finitary perspective, the existence of a Satisfaction Method as in (ii) above is also implicit in the usual expositions of Tarski's definitions.

\vspace{+1ex}
\textbf{A finitary perspective:} We highlight the word `\textit{\textbf{define}}' in (ii) above to emphasise the finitary perspective underlying this paper; which is that the concepts of `satisfaction' and `truth' under an interpretation are to be explicitly viewed as objective assignments by a convention that is witness-independent. A Platonist perspective would substitute `decide' for `define', thus implicitly suggesting that these concepts can `exist', in the sense of needing to be discovered by some witness-dependent means---eerily akin to a `revelation'---if the domain $\mathbb{D}$ is $N$.
\end{quote}

\noindent Classically, we can now inductively assign truth values of `satisfaction', `truth', and `falsity' to the compound formulas of  a first-order theory S under the interpretation $\mathcal{I}_{S(\mathbb{D})}$ in terms of \textit{only} the satisfiability of the atomic formulas of S over $\mathbb{D}$ as usual\footnote{See \cite{Me64}, p.51; \cite{Mu91}.}:

\begin{definition}
\label{sec:5.def.2}
A denumerable sequence $s$ of $\mathbb{D}$ satisfies $[\neg A]$ under $\mathcal{I}_{S(\mathbb{D})}$ if, and only if, $s$ does not satisfy $[A]$;
\end{definition}

\begin{definition}
\label{sec:5.def.3}
A denumerable sequence $s$ of $\mathbb{D}$ satisfies $[A \rightarrow B]$ under $\mathcal{I}_{S(\mathbb{D})}$ if, and only if, either it is not the case that $s$ satisfies $[A]$, or $s$ satisfies $[B]$;
\end{definition}

\begin{definition}
\label{sec:5.def.4}
A denumerable sequence $s$ of $\mathbb{D}$ satisfies $[(\forall x_{i})A]$ under $\mathcal{I}_{S(\mathbb{D})}$ if, and only if, given any denumerable sequence $t$ of $\mathbb{D}$ which differs from $s$ in at most the $i$'th component, $t$ satisfies $[A]$;
\end{definition}

\begin{quote}
\footnotesize We note that classical theory assumes without qualification\footnote{The need for a belief in the soundness of the Standard interpretation $\mathcal{I}_{PA(N,\ Standard)}$ of PA as an article of faith---falsified by Theorem \ref{sec:6.3.cor.3}---can be traced to this lack of qualification which, as Edward Nelson coloufully dramatises in \cite{Ne00}, ``\ldots is the battlefield where clash the armies of Platonists, intuitionists, and formalists. It differs from syntactical definitions because it invokes the notion of an infinite search. We have left the realm of the concrete for the speculative".} that, if the atomic formulas of $S$ are decidable, then it is always decidable whether a denumerable sequence $s$ of $D$ satisfies $[(\forall x_{i})A]$ under $\mathcal{I}_{S(D)}$.
\end{quote}

\begin{definition}
\label{sec:5.def.5}
A well-formed formula $[A]$ of $\mathbb{D}$ is true under $\mathcal{I}_{S(\mathbb{D})}$ if, and only if, given any denumerable sequence $t$ of $\mathbb{D}$, $t$ satisfies $[A]$;
\end{definition}

\begin{definition}
\label{sec:5.def.6}
A well-formed formula $[A]$ of $\mathbb{D}$ is false under $\mathcal{I}_{S(\mathbb{D})}$ if, and only if, it is not the case that $[A]$ is true under $\mathcal{I}_{S(\mathbb{D})}$.
\end{definition}

\noindent It follows that\footnote{cf. \cite{Me64}, pp.51-53.}:

\begin{theorem}
\label{sec:5.thm.1}
(\textit{Satisfaction Theorem}) If, for any interpretation $\mathcal{I}_{S(\mathbb{D})}$ of a first-order theory S, there is a Satisfaction Method SM($\mathcal{I}_{S(\mathbb{D})}$) which holds for a witness $\mathcal{W}_{\mathbb{D}}$ of $\mathbb{D}$, then:

\begin{quote}
(i) The $\Delta_{0}$ formulas of S are decidable as either true or false over $\mathbb{D}$ under $\mathcal{I}_{S(\mathbb{D})}$;

\vspace{+1ex}
(ii) If the $\Delta_{n}$ formulas of S are decidable as either true or as false over $\mathbb{D}$ under $\mathcal{I}_{S(\mathbb{D})}$, then so are the $\Delta(n+1)$ formulas of S.
\end{quote}
\end{theorem}

\vspace{+1ex}
\noindent \textbf{Proof} It follows from the above definitions that:

\vspace{+1ex}
\noindent (a) If, for any given atomic formula $[A(x_{1}, x_{2}, \ldots, x_{n})]$ of S, it is decidable by $\mathcal{W}_{\mathbb{D}}$ whether or not a given denumerable sequence $(a_{1}, a_{2}, \ldots)$ of $\mathbb{D}$ satisfies $[A(x_{1}, x_{2}, \ldots, x_{n})]$ in $\mathbb{D}$ under $\mathcal{I}_{S(\mathbb{D})}$ then, for any given compound formula $[A^{1}(x_{1}, x_{2}, \ldots, x_{n})]$ of S containing any one of the logical constants $\neg, \rightarrow, \forall$, it is decidable by $\mathcal{W}_{\mathbb{D}}$ whether or not $(a_{1}, a_{2}, \ldots)$ satisfies $[A^{1}(x_{1}, x_{2}, \ldots, x_{n})]$ in $\mathbb{D}$ under $\mathcal{I}_{S(\mathbb{D})}$;

\vspace{+1ex}
\noindent (b) If, for any given compound formula $[B^{n}(x_{1}, x_{2}, \ldots, x_{n})]$ of S containing $n$ of the logical constants $\neg, \rightarrow, \forall$, it is decidable by $\mathcal{W}_{\mathbb{D}}$ whether or not a given denumerable sequence $(a_{1}, a_{2}, \ldots)$ of $\mathbb{D}$ satisfies $[B^{n}(x_{1}, x_{2}, \ldots, x_{n})]$ in $\mathbb{D}$ under $\mathcal{I}_{S(\mathbb{D})}$ then, for any given compound formula $[B^{(n+1)}(x_{1}, x_{2}, \ldots, x_{n})]$ of S containing $n+1$ of the logical constants $\neg, \rightarrow, \forall$, it is decidable by $\mathcal{W}_{\mathbb{D}}$ whether or not $(a_{1}, a_{2}, \ldots)$ satisfies $[B^{(n+1)}(x_{1}, x_{2}, \ldots, x_{n})]$ in $\mathbb{D}$ under $\mathcal{I}_{S(\mathbb{D})}$;

\vspace{+1ex}
\noindent We thus have that:

\vspace{+1ex}
\noindent (c) The $\Delta_{0}$ formulas of S are decidable by $\mathcal{W}_{\mathbb{D}}$ as either true or false over $\mathbb{D}$ under $\mathcal{I}_{S(\mathbb{D})}$;

\vspace{+1ex}
\noindent (d) If the $\Delta_{n}$ formulas of S are decidable by $\mathcal{W}_{\mathbb{D}}$ as either true or as false over $\mathbb{D}$ under $\mathcal{I}_{S(\mathbb{D})}$, then so are the $\Delta(n+1)$ formulas of S. \hfill $\Box$

\vspace{+1ex}
\noindent In other words, if the atomic formulas of of S interpret under $\mathcal{I}_{S(\mathbb{D})}$ as decidable with respect to the Satisfaction Method SM($\mathcal{I}_{S(\mathbb{D})}$) by a witness $\mathcal{W}_{\mathbb{D}}$ over some domain $\mathbb{D}$, then the propositions of S (i.e., the $\Pi_{n}$ and $\Sigma_{n}$ formulas of S) also interpret as decidable with respect to SM($\mathcal{I}_{S(\mathbb{D})}$) by the witness $\mathcal{W}_{\mathbb{D}}$ over $\mathbb{D}$.

\vspace{+1ex}
\noindent We now consider the application of Tarski's definitions to various interpretations of first-order Peano Arithmetic PA.

\subsection{The standard interpretation of PA}
\label{sec:5.1}

The standard interpretation $\mathcal{I}_{PA(N,\ Standard)}$ of PA is obtained if, in $\mathcal{I}_{S(\mathbb{D})}$:

\begin{quote}
(a) we define S as PA with standard first-order predicate calculus as the underlying logic\footnote{Where the string $[(\exists \ldots)]$ is defined as---and is to be treated as an abbreviation for---the string $[\neg (\forall \ldots) \neg]$. We do not consider the case where the underlying logic is Hilbert's formalisation of Aristotle's logic of predicates in terms of his $\epsilon$-operator (\cite{Hi27}, pp.465-466).};

(b) we define $\mathbb{D}$ as $N$;

(c) for any atomic formula $[A(x_{1}, x_{2}, \ldots, x_{n})]$ of PA and sequence $(a_{1}, a_{2}, \ldots, a_{n})$ of $N$, we take $\|$SATCON($\mathcal{I}_{PA(N)}$)$\|$ as:

\begin{quote}
$\|$$A^{*}(a_{1}^{*}, a_{2}^{*}, \ldots, a_{n}^{*})$ holds in $N$ and, for any given sequence  $(b_{1}^{*}, b_{2}^{*}, \ldots, b_{n}^{*})$ of $N$, the proposition $A^{*}(b_{1}^{*}, b_{2}^{*}, \ldots, b_{n}^{*})$ is decidable in $N$$\|$;
\end{quote}

(d) we define the witness $\mathcal{W}_{(N,\ Standard)}$ informally as the `mathematical intuition' of a human intelligence for whom, classically, $\|$SATCON($\mathcal{I}_{PA(N)}$)$\|$ has been \textit{implicitly} accepted as \textit{objectively} `decidable' in $N$;

\begin{quote}
\footnotesize{
We shall show that such acceptance is justified, but needs to be made explicit since:

\begin{lemma}
\label{sec:5.1.lem.1}
$A^{*}(x_{1}, x_{2}, \ldots, x_{n})$ is both algorithmically verifiable and algorithmically computable in $N$ by $\mathcal{W}_{(N,\ Standard)}$. 
\end{lemma}

\vspace{+1ex}
\textbf{Proof} (i) It follows from the argument in Theorem \ref{sec:5.5.lem.1} (below) that $A^{*}(x_{1}, x_{2}, \ldots, x_{n})$ is algorithmically verifiable in $N$ by $\mathcal{W}_{(N,\ Standard)}$.

\vspace{+.5ex}
(ii) It follows from the argument in Theorem \ref{sec:5.5.lem.2} (below) that $A^{*}(x_{1}, x_{2},$ $\ldots, x_{n})$ is algorithmically computable in $N$ by $\mathcal{W}_{(N,\ Standard)}$. The lemma follows.\hfill $\Box$

\vspace{+1ex}
Now, although it is not immediately obvious from the standard interpretation of PA which of (i) or (ii) may be taken for \textit{explicitly} deciding $\|$SATCON($\mathcal{I}_{PA(N)}$)$\|$ by the witness $\mathcal{W}_{(N,\ Standard)}$, we shall show in Section \ref{sec:5.3} that (i) is consistent with (e) below; and in Section \ref{sec:5.5} that (ii) is inconsistent with (e). Thus the standard interpretation of PA implicitly presumes (i).
}
\end{quote}

(e) we postulate that Aristotle's particularisation holds over $N$\footnote{Hence a PA formula such as $[(\exists x)F(x)]$ interprets under $\mathcal{I}_{PA(N,\ Standard)}$ as `There is some natural number $n$ such that $F(n)$ holds in $N$.}.
\end{quote}

\vspace{+1ex}
\noindent Clearly, (e) does not form any part of Tarski's inductive definitions of the satisfaction, and truth, of the formulas of PA under the above interpretation. Moreover, its inclusion makes $\mathcal{I}_{PA(N,\ Standard)}$ extraneously non-finitary\footnote{\cite{Br08}.}.

\begin{quote}
\footnotesize{
The question arises: Can we formulate the `standard' interpretation of PA without assuming (e) extraneously?

We answer this question affirmatively in Section \ref{sec:5.3} where:

\vspace{+.5ex}
(1) We replace the `mathematical intuition' of a human intelligence by defining an `objective' witness $\mathcal{W}_{(\mathbb{N},\ Instantiational)}$ as the meta-theory $\mathcal{M}_{PA}$ of PA;

\vspace{+.5ex}
(2) We show that $\mathcal{W}_{(\mathbb{N},\ Instantiational)}$ can decide whether $c^{*}=d^{*}$ is true or false by instantiationally computing the Boolean function $A^{*}(x_{1}, x_{2}, \ldots, x_{n})$ for any given sequence of natural numbers $(b^{*}_{1}, b^{*}_{2}, \ldots, b^{*}_{n})$;

\vspace{+.5ex}
(3) We show that this yields an instantiational interpretation of PA over the structure $[\mathbb{N}]$ of the PA numerals that is sound if, and only if, (e) holds.

\vspace{+.5ex}
(4) $\mathcal{W}_{(\mathbb{N},\ Instantiational)}$ is thus an instantiational formulation of the standard interpretation of PA over $\mathbb{N}$ (which is presumed to be sound).
}
\end{quote}

\noindent We note further that if PA is $\omega$-\textit{in}consistent, then Aristotle's particularisation does not hold over $N$, and the interpretation $\mathcal{I}_{PA(N,\ Standard)}$ is not sound.

\subsection{G\"{o}del's non-standard interpretation of PA}
\label{sec:5.2}

A non-standard (G\"{o}delian) interpretation $\mathcal{I}_{PA(N_{\omega},\ Non-standard)}$ of a putative $\omega$-consistent PA is obtained if, in $\mathcal{I}_{S(\mathbb{D})}$:

\begin{quote}
(a) we define S as PA with standard first-order predicate calculus as the underlying logic;

(b) we define $\mathbb{D}$ as an undefined extension $N_{\omega}$ of $N$;

(c) for any atomic formula $[A(x_{1}, x_{2}, \ldots, x_{n})]$ of PA and sequence $(a_{1}, a_{2}, \ldots, a_{n})$ of $N_{\omega}$, we take $\|$SATCON($\mathcal{I}_{PA(N_{\omega})}$)$\|$ as:

\begin{quote}
$\|$$A^{*}(a_{1}^{*}, a_{2}^{*}, \ldots, a_{n}^{*})$ holds in $N_{\omega}$ and, for any given sequence  $(b_{1}^{*}, b_{2}^{*}, \ldots, b_{n}^{*})$ of $N_{\omega}$, the proposition $A^{*}(b_{1}^{*}, b_{2}^{*}, \ldots, b_{n}^{*})$ is decidable as either holding or not holding in $N_{\omega}$$\|$;
\end{quote}

(d) we \textit{postulate} that $\|$SATCON($\mathcal{I}_{PA(N_{\omega})}$)$\|$ is always decidable by a putative witness $\mathcal{W}_{\mathbb{N_{\omega}}}$, and that $\mathcal{W}_{\mathbb{N_{\omega}}}$ can, further, determine some numbers in $N_{\omega}$ which are not natural numbers;

(e) we assume that PA is $\omega$-consistent.
\end{quote}

\noindent Clearly, the interpretation $\mathcal{I}_{PA(N_{\omega},\ Non-standard)}$ of a putative $\omega$-consistent PA cannot claim to be finitary. Moreover, if PA is $\omega$-\textit{in}consistent, then the G\"{o}delian non-standard interpretation $\mathcal{I}_{PA(N_{\omega},\ Non-standard)}$ of PA is also not sound\footnote{In which case we cannot validly conclude from G\"{o}del's formal reasoning in (\cite{Go31}) that PA must have a non-standard model.}.

\subsection{An instantiational interpretation of PA in PA}
\label{sec:5.3}

We next consider the instantiational interpretation\footnote{The raison d'\^{e}tre, and significance, of such interpretation is outlined in this short \href{http://alixcomsi.com/8_Meeting_Wittgenstein_requirement_1000.pdf}{note}.} $\mathcal{I}_{PA(\mathbb{N},\ Instantiational)}$ of PA where:

\begin{quote}
(a) we define S as PA with standard first-order predicate calculus as the underlying logic;

(b) we define $D$ as the set $\mathbb {N}$ of PA numerals;

(c) for any atomic formula $[A(x_{1}, x_{2}, \ldots, x_{n})]$ of PA and any sequence $[(a_{1}, a_{2}, \ldots, a_{n})]$ of PA numerals, we take $\|$SATCON($\mathcal{I}_{PA(\mathcal{PA})}$)$\|$ as:

\begin{quote}
$\|$$[A(a_{1}, a_{2}, \ldots, a_{n})]$ is provable in PA and, for any given sequence of numerals $[(b_{1}, b_{2}, \ldots, b_{n})]$ of PA, the formula $[A(b_{1}, b_{2}, \ldots, b_{n})]$ is decidable as either provable or not provable in PA$\|$;
\end{quote}

(d) we define the witness $\mathcal{W}_{(\mathbb{N},\ Instantiational)}$ as the meta-theory $\mathcal{M}_{PA}$ of PA.

\begin{quote}
\footnotesize{
\begin{lemma}
\label{sec:5.3.lem.1}
$[A(x_{1}, x_{2}, \ldots, x_{n})]$ is always algorithmically verifiable in PA by $\mathcal{W}_{(\mathbb{N},\ Instantiational)}$. 
\end{lemma}

\vspace{+1ex}
\textbf{Proof} It follows from G\"{o}del's definition of the primitive recursive relation $xBy$\footnote{\cite{Go31}, p. 22(45).}---where $x$ is the G\"{o}del number of a proof sequence in PA whose last term is the PA formula with G\"{o}del-number $y$---that, if $[A(x_{1}, x_{2}, \ldots, x_{n})]$ is an atomic formula of PA, then $\mathcal{M}_{PA}$ can algorithmically verify for any given sequence $[(b_{1}, b_{2}, \ldots, b_{n})]$ of PA numerals which one of the PA formulas $[A(b_{1}, b_{2}, \ldots, b_{n})]$ and $[\neg A(b_{1}, b_{2}, \ldots, b_{n})]$ is necessarily PA-provable.\hfill $\Box$
}
\end{quote}
\end{quote}

\noindent Now, if PA is consistent but not $\omega$-consistent, then there is a G\"{o}delian formula $[R(x)]$ such that (see Section \ref{sec:6.3.1}):

\begin{quote}
(i) $[(\forall x)R(x)]$ is not PA-provable;

(ii) $[\neg (\forall x)R(x)]$ is PA-provable;

(iii) for any given numeral $[n]$, $[R(n)]$ is PA-provable.
\end{quote}

\noindent However, if $\mathcal{I}_{PA(\mathbb{N},\ Instantiational)}$ is sound, then (ii) implies contradictorily that it is not the case that, for any given numeral $[n]$, $[R(n)]$ is PA-provable.

\vspace{+1ex}
\noindent It follows that if $\mathcal{I}_{PA(\mathbb{N},\ Instantiational)}$ is sound then PA is $\omega$-consistent and, ipso facto, Aristotle's particularisation must hold over $\mathbb{N}$.

\vspace{+1ex}
\noindent Moreover, if PA is consistent, then every PA-provable formula interprets as true under some sound interpretation of PA. Hence $\mathcal{M}_{PA}$ can effectively decide whether, for any given sequence of natural numbers $(b_{1}^{*}, b_{2}^{*},$ $\ldots, b_{n}^{*})$ in $\mathbb{N}$, the proposition $A^{*}(b_{1}^{*}, b_{2}^{*}, \ldots, b_{n}^{*})$ holds or not in $\mathbb{N}$.

\vspace{+1ex}
\noindent It follows that $\mathcal{I}_{PA(\mathbb{N},\ Instantiational)}$ is an instantiational formulation of the `standard' interpretation of PA in which we do not need to extraneously assume that Aristotle's particularisation holds over $\mathbb{N}$.

\vspace{+1ex}
\noindent The interpretation $\mathcal{I}_{PA(\mathbb{N},\ Instantiational)}$ is of interest because, if it were a sound interpretation of PA, then PA would meta-mathematically establish its own consistency\footnote{This is not possible mathematically by G\"{o}del's Theorem XI in \cite{Go31}, p.36, if PA is $\omega$-consistent.}!

\subsection{A set-theoretic interpretation of PA}
\label{sec:5.4}

We consider next a set-theoretic interpretation $\mathcal{I}_{PA(\mathcal{ZF},\ Cantor)}$ of PA over the domain of ZF sets, which is obtained if, in $\mathcal{I}_{S(\mathbb{D})}$:

\begin{quote}
(a) we define S as PA with standard first-order predicate calculus as the underlying logic;

(b) we define $\mathbb{D}$ as ZF;

(c) for any atomic formula $[A(x_{1}, x_{2}, \ldots, x_{n})]$ and sequence $[(a_{1}, a_{2}, \ldots, a_{n})]$ of PA, we take $\|$SATCON($\mathcal{I}_{S(\mathcal{ZF})}$)$\|$ as:

\begin{quote}
$\|$$[A^{*}(a_{1}^{*}, a_{2}^{*}, \ldots, a_{n}^{*})]$ is provable in ZF and, for any given sequence  $[(b_{1}^{*}, b_{2}^{*}, \ldots, b_{n}^{*})]$ of ZF, the formula $[A^{*}(b_{1}^{*}, b_{2}^{*}, \ldots, b_{n}^{*})]$ is decidable as either provable or not provable in ZF$\|$;
\end{quote}

(d) we define the witness $\mathcal{W}_{\mathcal{ZF}}$ as the meta-theory $\mathcal{M}_{ZF}$ of ZF which can always decide effectively whether or not $\|$SATCON($\mathcal{I}_{S(\mathcal{ZF})}$)$\|$ holds in ZF.
\end{quote}

\noindent Now, if the set-theoretic interpretation $\mathcal{I}_{PA(\mathcal{ZF},\ Cantor)}$ of PA is sound, then every sound interpretation of ZF would, ipso facto, be a sound interpretation of PA. In Appendix E, Section \ref{sec:10.1} I show, however, that this is not the case, and so the set-theoretic interpretation $\mathcal{I}_{PA(\mathcal{ZF},\ Cantor)}$ of PA is not sound.

\subsection{A purely algorithmic interpretation of PA}
\label{sec:5.5}

We finally consider the purely algorithmic interpretation $\mathcal{I}_{PA(N,\ Algorithmic)}$ of PA where:

\begin{quote}
(a) we define S as PA with standard first-order predicate calculus as the underlying logic;

(b) we define $\mathbb{D}$ as $N$;

(c) for any atomic formula $[A(x_{1}, x_{2}, \ldots, x_{n})]$ and sequence $(a_{1}, a_{2}, \ldots, a_{n})$ of natural numbers in $N$, we take $\|$SATCON($\mathcal{I}_{PA(N)}$)$\|$ as:

\begin{quote}
$\|$$A^{*}(a_{1}^{*}, a_{2}^{*}, \ldots, a_{n}^{*})$ holds in $N$ and, for any given sequence  $(b_{1}^{*}, b_{2}^{*}, \ldots, b_{n}^{*})$ of $N$, the proposition $A^{*}(b_{1}^{*}, b_{2}^{*}, \ldots, b_{n}^{*})$ is decidable as either holding or not holding in $N$$\|$;
\end{quote}

(d) we define the witness $\mathcal{W}_{(N,\ Algorithmic)}$ as any simple functional language that computes $[A(x_{1}, x_{2}, \ldots, x_{n})]$ and gives evidence that $\|$SATCON($\mathcal{I}_{PA(N)}$)$\|$ is always \textit{effectively} decidable in $N$:

\begin{quote}
\footnotesize
\begin{lemma}
\label{sec:5.5.lem.1}
$A^{*}(x_{1}, x_{2}, \ldots, x_{n})$ is always algorithmically computable in $N$ by $\mathcal{W}_{(N,\ Algorithmic)}$. 
\end{lemma}

\vspace{+1ex}
\textbf{Proof} If $[A(x_{1}, x_{2}, \ldots, x_{n})]$ is an atomic formula of PA then, for any given sequence of numerals $[b_{1}, b_{2}, \ldots, b_{n}]$, the PA formula $[A(b_{1}, b_{2},$ $\ldots, b_{n})]$ is an atomic formula of the form $[c=d]$, where $[c]$ and $[d]$ are atomic PA formulas that denote PA numerals. Since $[c]$ and $[d]$ are recursively defined formulas in the language of PA, it follows from a standard result\footnote{For any natural numbers $m,\ n$, if $m \neq n$, then PA proves $[\neg(m = n)]$ (\cite{Me64}, p.110, Proposition 3.6). The converse is obviously true.} that, if PA is consistent, then $[c=d]$ is algorithmically computable as either true or false in $N$. In other words, if PA is consistent, then $[A(x_{1}, x_{2}, \ldots, x_{n})]$ is algorithmically computable (since there is a deterministic algorithm that, for any given sequence of numerals $[b_{1}, b_{2}, \ldots, b_{n}]$, will give evidence whether $[A(b_{1}, b_{2},$ $\ldots, b_{n})]$ interprets as true or false in $N$. The lemma follows.\hfill $\Box$
\end{quote}
\end{quote}

\noindent It follows that $\mathcal{I}_{PA(N,\ Algorithmic)}$ is an algorithmic formulation of the `standard' interpretation of PA in which we do not extraneously assume that Aristotle's particularisation holds over $N$.

\begin{quote}
\footnotesize{
We shall show that if $\mathcal{I}_{PA(N,\ Algorithmic)}$ is sound, then PA is \textit{not} $\omega$-consistent. Hence Aristotle's particularisation does not hold over $N$, and the interpretation is finitary and intuitionistically unobjectionable. Moreover---since the Law of the Excluded Middle is provable in an $\omega$-inconsistent PA (and therefore holds in $N$)---it achieves this without the discomforting, stringent, Intuitionistic requirement that we reject the underlying logic of PA!
}
\end{quote}

\section{Formally defining the standard interpretation of PA \textit{finitarily}}
\label{sec:5.4.0.0}

It follows from the analysis of the classical applicability of Tarski's inductive definitions of `satisfiability' and `truth' in Section \ref{sec:5.4.0} that we can formally define---as detailed in Section \ref{sec:7}, Appendix B---the standard interpretation $\mathcal{I}_{PA(N,\ Standard)}$ of PA where:

\begin{quote}
(a) we define S as PA with standard first-order predicate calculus as the underlying logic;

\vspace{+1ex}
(b) we define $\mathbb{D}$ as $N$;

\vspace{+1ex}
(c) we take SM($\mathcal{I}_{PA(N,\ Standard)}$) as any simple functional language.
\end{quote}

\noindent We note that:

\begin{theorem}
\label{sec:5.5.lem.1}
The atomic formulas of PA are algorithmically verifiable under the standard interpretation $\mathcal{I}_{PA(N,\ Standard)}$. 
\end{theorem}

\noindent \textbf{Proof} If $[A(x_{1}, x_{2}, \ldots, x_{n})]$ is an atomic formula of PA then, for any given denumerable sequence of numerals $[b_{1}, b_{2}, \ldots]$, the PA formula $[A(b_{1}, b_{2},$ $\ldots, b_{n})]$ is an atomic formula of the form $[c=d]$, where $[c]$ and $[d]$ are atomic PA formulas that denote PA numerals. Since $[c]$ and $[d]$ are recursively defined formulas in the language of PA, it follows from a standard result that, if PA is consistent, then $[c=d]$ interprets as the proposition $c=d$ which either holds or not for a witness $\mathcal{W}_{N}$ in $N$.

\vspace{+1ex}
\noindent Hence, if PA is consistent, then $[A(x_{1}, x_{2}, \ldots, x_{n})]$ is algorithmically verifiable since, for any given denumerable sequence of numerals $[b_{1}, b_{2}, \ldots]$, we can define a deterministic algorithm that provides evidence that the PA formula $[A(b_{1}, b_{2}, \ldots, b_{n})]$ is decidable under the interpretation.

\vspace{+1ex}
\noindent The theorem follows.\hfill $\Box$

\vspace{+1ex}
\noindent It immediately follows that:

\begin{corollary}
\label{sec:5.5.cor.1}
The `satisfaction' and `truth' of PA formulas containing logical constants can be defined under the standard interpretation of PA in terms of the evidence provided by the computations of a simple functional language.
\end{corollary}

\begin{corollary}
\label{sec:5.5.cor.2}
The PA-formulas are decidable under the standard interpretation of PA if, and only if, they are algorithmically verifiable under the interpretation.
\end{corollary}

\subsection{Defining `algorithmic truth' under the standard interpretation of PA}
\label{sec:6.2}

Now we note that, in addition to Theorem \ref{sec:5.5.lem.1}:

\begin{theorem}
\label{sec:5.5.lem.2}
The atomic formulas of PA are algorithmically computable under the standard interpretation $\mathcal{I}_{PA(N,\ Standard)}$. 
\end{theorem}

\noindent \textbf{Proof} If $[A(x_{1}, x_{2}, \ldots, x_{n})]$ is an atomic formula of PA then we can define a deterministic algorithm that, for any given denumerable sequence of numerals $[b_{1}, b_{2}, \ldots]$, provides evidence whether the PA formula $[A(b_{1}, b_{2}, \ldots, b_{n})]$ is true or false under the interpretation.

\vspace{+1ex}
\noindent The theorem follows.\hfill $\Box$

\vspace{+1ex}
\noindent This suggests the following definitions:

\begin{definition}
\label{sec:6.2.def.1}
A well-formed formula $[A]$ of PA is algorithmically true under $\mathcal{I}_{PA(N,\ Standard)}$ if, and only if, there is a deterministic algorithm which provides evidence that, given any denumerable sequence $t$ of $N$, $t$ satisfies $[A]$;
\end{definition}

\begin{definition}
\label{sec:6.2.def.2}
A well-formed formula $[A]$ of PA is algorithmically false under $\mathcal{I}_{PA(N,\ Standard)}$ if, and only if, it is not algorithmically true under $\mathcal{I}_{PA(N)}$.
\end{definition}

\subsubsection{An algorithmic interpretation of the PA axioms}
\label{sec:6.2.1}

The significance of defining `algorithmic truth' under $\mathcal{I}_{PA(N,\ Standard)}$ as above is that:

\begin{lemma}
\label{sec:6.2.lem.1}
The PA axioms PA$_{1}$ to PA$_{8}$ are algorithmically computable as algorithmically true over $N$ under the interpretation $\mathcal{I}_{PA(N,\ Standard)}$.
\end{lemma}

\vspace{+1ex}
\noindent \textbf{Proof} Since $[x+y]$, $[x \star y]$, $[x = y]$, $[{x^{\prime}}]$ are defined recursively\footnote{cf. \cite{Go31}, p.17.}, the PA axioms PA$_{1}$ to PA$_{8}$ interpret as recursive relations that do not involve any quantification. The lemma follows straightforwardly from Definitions \ref{sec:5.def.1} to \ref{sec:5.def.6} in Section \ref{sec:5.4.0} and Theorem \ref{sec:5.5.lem.1}.\hfill $\Box$

\begin{lemma}
\label{sec:6.2.lem.2}
For any given PA formula $[F(x)]$, the Induction axiom schema $[F(0)$ $\rightarrow (((\forall x)(F(x) \rightarrow F(x^{\prime}))) \rightarrow (\forall x)F(x))]$ interprets as algorithmically true under $\mathcal{I}_{PA(N,\ Standard)}$.
\end{lemma}

\vspace{+1ex}
\noindent \textbf{Proof} By Definitions \ref{sec:5.def.1} to \ref{sec:6.2.def.2}:

\begin{quote}

(a) If $[F(0)]$ interprets as algorithmically false under $\mathcal{I}_{PA(N,\ Standard)}$ the lemma is proved.

\begin{quote}
\footnotesize Since $[F(0) \rightarrow (((\forall x)(F(x) \rightarrow F(x^{\prime}))) \rightarrow (\forall x)F(x))]$ interprets as algorithmically true if, and only if, either $[F(0)]$ interprets as algorithmically false or $[((\forall x)(F(x) \rightarrow F(x^{\prime}))) \rightarrow (\forall x)F(x)]$ interprets as algorithmically true.
\end{quote}

\vspace{+1ex}
(b) If $[F(0)]$ interprets as algorithmically true and $[(\forall x)(F(x) \rightarrow F(x^{\prime}))]$ interprets as algorithmically false under $\mathcal{I}_{PA(N,\ Standard)}$, the lemma is proved.

\vspace{+1ex}
(c) If $[F(0)]$ and $[(\forall x)(F(x) \rightarrow F(x^{\prime}))]$ both interpret as algorithmically true under $\mathcal{I}_{PA(N,\ Standard)}$, then by Definition \ref{sec:6.2.def.1} there is a deterministic Turing machine that computes $[F(x)]$ and, for any natural number $n$, will give evidence that the formula $[F(n) \rightarrow F(n^{\prime})]$ is true under $\mathcal{I}_{PA(N,\ Standard)}$.

\vspace{+1ex}
Since $[F(0)]$ interprets as algorithmically true under $\mathcal{I}_{PA(N,\ Standard)}$, it follows that there is a deterministic Turing machine that computes $[F(x)]$ and, for any natural number $n$, will give evidence that the formula $[F(n)]$ is true under the interpretation.

\vspace{+1ex}
Hence $[(\forall x)F(x)]$ is algorithmically true under $\mathcal{I}_{PA(N,\ Standard)}$.
\end{quote}

\noindent Since the above cases are exhaustive, the lemma follows.\hfill $\Box$

\begin{quote}
\footnotesize
\textbf{The Poincar\'{e}-Hilbert debate:} We note that Lemma \ref{sec:6.2.lem.2} appears to settle the Poincar\'{e}-Hilbert debate\footnote{See \cite{Hi27}, p.472; also \cite{Br13}, p.59; \cite{We27}, p.482; \cite{Pa71}, p.502-503.} in the latter's favour. Poincar\'{e} believed that the Induction Axiom could not be justified finitarily, as any such argument would necessarily need to appeal to infinite induction. Hilbert believed that a finitary proof of the consistency of PA was possible. 
\end{quote}

\begin{lemma}
\label{sec:6.2.lem.3}
Generalisation preserves algorithmic truth under $\mathcal{I}_{PA(N,\ Standard)}$.
\end{lemma}

\noindent \textbf{Proof} The two meta-assertions:

\begin{quote}
`$[F(x)]$ interprets as algorithmically true under $\mathcal{I}_{PA(N,\ Standard)}$\footnote{See Definition \ref{sec:5.def.5}}'

\vspace{+.5ex}
and

\vspace{+.5ex}
`$[(\forall x)F(x)]$ interprets as algorithmically true under $\mathcal{I}_{PA(N,\ Standard)}$'
\end{quote}

\noindent both mean:

\begin{quote}
$[F(x)]$ is algorithmically computable as always true under $\mathcal{I}_{PA(N),}$ $_{Standard)}$. \hfill $\Box$
\end{quote} 

\vspace{+1ex}
\noindent It is also straightforward to see that:

\begin{lemma}
\label{sec:6.2.lem.4}
Modus Ponens preserves algorithmic truth under $\mathcal{I}_{PA(N,\ Standard)}$. \hfill $\Box$
\end{lemma}

\noindent We thus have that:

\begin{theorem}
\label{sec:6.2.lem.5}
The axioms of PA are always algorithmically true under the interpretation $\mathcal{I}_{PA(N,\ Standard)}$, and the rules of inference of PA preserve the properties of algorithmic satisfaction/truth under $\mathcal{I}_{PA(N,\ Standard)}$\footnote{Without appeal, moreover, to Aristotle's particularisation.}.\hfill $\Box$
\end{theorem}

\subsubsection{The algorithmic interpretation $\mathcal{I}_{PA(N,\ Algorithmic)}$ of PA over $N$ is sound}
\label{sec:6.2.2}

We conclude that there is a deterministic algorithmic interpretation $\mathcal{I}_{PA(N,\ Algorithmic)}$ of PA over $N$---formally defined in Section \ref{sec:8}, Appendix C---such that:

\begin{theorem}
\label{sec:6.2.thm.1}
The interpretation $\mathcal{I}_{PA(N,\ Algorithmic)}$ of PA is sound\footnote{In the sense of Definitions \ref{sec:A.def.9} and \ref{sec:A.def.10}.}.
\end{theorem}

\noindent \textbf{Proof} It follows immediately from Theorem \ref{sec:6.2.lem.5} and Section \ref{sec:8}, Appendix C, that the axioms of PA are always true under the interpretation $\mathcal{I}_{PA(N,\ Algorithmic)}$, and the rules of inference of PA preserve the properties of satisfaction/truth under $\mathcal{I}_{PA(N,\ Algorithmic)}$.\hfill $\Box$

\vspace{+1ex}
\noindent We thus have a finitary proof that:

\begin{theorem}
\label{sec:6.2.thm.2}
PA is consistent. \hfill $\Box$
\end{theorem}

\begin{quote}
\footnotesize
\textbf{Hilbert's Second Problem:} We note---but do not consider further as it is not germane to the intent of this investigation---that Lemma \ref{sec:6.2.thm.2} offers a partial resolution to Hilbert's Second Problem, which asks for a finitary proof that the second order Arithmetical axioms are consistent\footnote{``When we are engaged in investigating the foundations of a science, we must set up a system of axioms which contains an exact and complete description of the relations subsisting between the elementary ideas of that science. \ldots But above all I wish to designate the following as the most important among the numerous questions which can be asked with regard to the axioms: To prove that they are not contradictory, that is, that a definite number of logical steps based upon them can never lead to contradictory results. In geometry, the proof of the compatibility of the axioms can be effected by constructing a suitable field of numbers, such that analogous relations between the numbers of this field correspond to the geometrical axioms. \ldots On the other hand a direct method is needed for the proof of the compatibility of the arithmetical axioms." \ldots \cite{Nw02}.}. 
\end{quote}

\section{A Provability Theorem for PA}
\label{sec:6.3}

We now show that PA can have no non-standard model\footnote{We consider the usual arguments for the existence of non-standard models of PA in Section \ref{sec:11}, Appendix F.}, since it is `algorithmically' complete in the sense that:

\begin{theorem}
\label{sec:6.3.thm.1}
(Provability Theorem for PA) A PA formula $[F(x)]$ is PA-provable if, and only if, $[F(x)]$ is algorithmically computable as always true in $N$.
\end{theorem}

\vspace{+1ex}
\noindent \textbf{Proof} We have by definition that $[(\forall x)F(x)]$ interprets as true under the interpretation $\mathcal{I}_{PA(N,\ Algorithmic)}$ if, and only if, $[F(x)]$ is algorithmically computable as always true in $N$.

\vspace{+1ex}
\noindent Since $\mathcal{I}_{PA(N,\ Algorithmic)}$ is sound, it defines a finitary model of PA over $N$---say $\mathcal{M}_{PA(\beta)}$---such that:

\begin{itemize}
\item
If $[(\forall x)F(x)]$ is PA-provable, then $[F(x)]$ is algorithmically computable as always true in $N$;

\vspace{+1ex}
\item
\noindent If $[\neg(\forall x)F(x)]$ is PA-provable, then it is not the case that $[F(x)]$ is algorithmically computable as always true in $N$.
\end{itemize}

\noindent Now, we cannot have that both $[(\forall x)F(x)]$ and $[\neg(\forall x)F(x)]$ are PA-unprovable for some PA formula $[F(x)]$, as this would yield the contradiction:

\begin{itemize}
\item
There is a finitary model---say $M1_{\beta}$---of PA+$[(\forall x)F(x)]$ in which $[F(x)]$ is algorithmically computable as always true in $N$; and

\vspace{+1ex}
\item
There is a finitary model---say $M2_{\beta}$---of PA+$[\neg(\forall x)F(x)]$ in which it is not the case that $[F(x)]$ is algorithmically computable as always true in $N$.
\end{itemize}

\noindent Further, we cannot have that:

\begin{itemize}
\item
$[F(x)]$ is algorithmically computable as always true in $N$, and $[\neg(\forall x)F(x)]$ is PA-provable;
\end{itemize}

\noindent nor that:

\begin{itemize}
\item
It is not the case that $[F(x)]$ is algorithmically computable as always true in $N$, and $[(\forall x)F(x)]$ is PA-provable.
\end{itemize}

\noindent The lemma follows.\hfill $\Box$

\vspace{+1ex}
\noindent We conclude that:

\begin{corollary}
\label{sec:6.3.thm.1.cor.1}
The provable formulas of PA are precisely those that are algorithmically computable as always true under a sound interpretation of PA.
\end{corollary}

\vspace{+1ex}
\noindent We further conclude that\footnote{cf. Hilbert's remarks at the International Congress of Mathematicians at Paris in 1900: ``The axioms of arithmetic are essentially nothing else than the known rules of calculation, with the addition of the axiom of continuity. I recently collected them and in so doing replaced the axiom of continuity by two simpler axioms, namely, the well-known axiom of Archimedes, and a new axiom essentially as follows: that numbers form a system of things which is capable of no further extension, as long as all the other axioms hold (axiom of completeness)." \ldots \cite{Nw02}.}:

\begin{corollary}
\label{sec:6.3.thm.1.cor.1}
PA is categorical.
\end{corollary}

\section{PA is \textit{not} $\omega$-consistent}
\label{sec:6.3.1}

In his seminal 1931 paper on formally undecidable arithmetical propositions\footnote{\cite{Go31}.}, G\"{o}del showed that\footnote{\cite{Go31}, Theorem VI, p.24.}:

\begin{lemma}
\label{sec:6.3.1.lem.1}
If a Peano Arithmetic such as PA is $\omega$-consistent, then there is a constructively definable PA-formula $[R(x)]$\footnote{In his argument, G\"{o}del refers to this formula only by its `G\"{o}del' number `\(r\)'; \cite{Go31}, p.25, Eqn.(12).} such that neither $[(\forall x)R(x)$) nor $[\neg(\forall x)R(x)]$ are PA-provable\footnote{\cite{Go31}, p.25(1) \& p.26(2).}.$\Box$
\end{lemma}

\noindent G\"{o}del concluded that:

\begin{lemma}
\label{sec:6.3.1.lem.2}
Any $\omega\)-consistent Peano Arithmetic such as PA has a consistent, but $\omega\)-\textit{in}consistent, extension PA$'$, obtained by adding the formula $[\neg(\forall x)R(x)]$ as an axiom to PA\footnote{\cite{Go31}, p.27.}.$\Box$
\end{lemma}

\noindent Specifically, G\"{o}del's reasoning shows that:

\begin{lemma}
\label{sec:6.3.1.lem.3}
If PA is consistent and $[(\forall x)R(x)]$ is assumed PA-provable, then $[\neg(\forall x)R(x)]$ is PA-provable\footnote{This follows from G\"{o}del's argument in \cite{Go31}, p.26(1).}.$\Box$
\end{lemma}

\begin{lemma}
\label{sec:6.3.1.lem.4}
If PA is $\omega$-consistent and $[\neg(\forall x)R(x)]$ is assumed PA-provable, then $[(\forall x)R(x)]$ is PA-provable\footnote{This follows from G\"{o}del's argument in \cite{Go31}, p.26(2).}.$\Box$
\end{lemma}

\vspace{+1ex}
\noindent However, by the argument in Theorem \ref{sec:6.3.thm.1} it now follows that:

\begin{corollary}
\label{sec:6.3.cor.1}
The PA formula $[\neg(\forall x)R(x)]$ is PA-provable. \hfill $\Box$
\end{corollary}

\begin{quote}
\footnotesize Of course $[\neg(\forall x)R(x)]$ interprets under $\mathcal{I}_{PA(N,\ Algorithmic)}$ as the assertion:

\begin{quote}
There is no deterministic algorithm that will compute $[R(x)]$ and, for any given natural number $n$, provide evidence that $R^{*}(n)$ is a true arithmetical proposition in $N$.
\end{quote}

However, since G\"{o}del has shown that the PA-formula $[R(n)]$ is PA-provable for any given PA-numeral $[n]$, it follows that:

\begin{quote}
For any given natural number $n$, there is always some deterministic algorithm that will compute $[R(n)]$ and provide evidence that $R^{*}(n)$ is a true arithmetical proposition in $N$.
\end{quote}

Thus the PA-formula $[(\forall x)R(x)]$ is algorithmically verifiable as true over $N$, but not algorithmically computable as true over $N$. The arithmetical relation $R^{*}(x)$ is thus a Halting-type of relation, such that although $R^{*}(x)$ is a tautology over $N$, there is no deterministic algorithm that will compute $[R(x)]$ and, for any given natural number $n$, give evidence that $R^{*}(n)$ is a true arithmetical proposition in $N$.
\end{quote}

\vspace{+1ex}
\noindent We conclude that:

\begin{corollary}
\label{sec:6.3.cor.2}
PA is \textit{not} $\omega$-consistent.\footnote{This conclusion is contrary to accepted dogma. See, for instance, Davis' remarks in \cite{Da82}, p.129(iii) that ``\ldots there is no equivocation. Either an adequate arithmetical logic is $\omega$-inconsistent (in which case it is possible to prove false statements within it) or it has an unsolvable decision problem and is subject to the limitations of G\"{o}del's incompleteness theorem".}
\end{corollary}

\vspace{+1ex}
\noindent \textbf{Proof} G\"{o}del has shown that if PA is consistent, then $[R(n)]$ is PA-provable for any given PA numeral $[n]$\footnote{\cite{Go31}, p.26(2).}. By Corollary \ref{sec:6.3.cor.1} and the definition of $\omega$-consistency, if PA is consistent then it is \textit{not} $\omega$-consistent.\hfill $\Box$

\begin{corollary}
\label{sec:6.3.cor.3}
The standard interpretation $\mathcal{I}_{PA(N,\ Standard)}$ of PA is not sound\footnote{In the sense of Definitions \ref{sec:A.def.9} and \ref{sec:A.def.10}.}, and does not yield a model of PA\footnote{We note that finitists of all hues---ranging from Brouwer \cite{Br08} to Alexander Yessenin-Volpin \cite{He04}---have persistently questioned the soundness of the `standard' interpretation $\mathcal{I}_{PA(N,\ Standard)}$.}.
\end{corollary}

\vspace{+1ex}
\noindent \textbf{Proof} By Corollary \ref{sec:4.2.cor.1} if PA is consistent but not $\omega$-consistent, then Aristotle's particularisation does not hold over $N$. Since the `standard', interpretation of PA appeals to Aristotle's particularisation, the lemma follows.\hfill $\Box$

\section{Appendix A: The significance of $\omega$-consistency and Hilbert's program}
\label{sec:4}

In order to avoid intuitionistic objections to his reasoning in his seminal 1931 paper on formally undecidable arithmetical propositions\footnote{\cite{Go31}.}, Kurt G\"{o}del did not assume that the standard interpretation $\mathcal{I}_{PA(N,\ Standard)}$ of PA is sound\footnote{In the sense of Definitions \ref{sec:A.def.9} and \ref{sec:A.def.10}.}. Instead, G\"{o}del introduced the syntactic property of $\omega$-consistency\footnote{Definition \ref{sec:A.def.8}.} as an explicit assumption in his formal reasoning\footnote{\cite{Go31}, p.23 and p.28.}. G\"{o}del explained at some length\footnote{In his introduction on p.9 of \cite{Go31}.} that his reasons for introducing $\omega$-consistency as an explicit assumption in his formal reasoning was to avoid appealing to the semantic concept of classical arithmetical truth---a concept which is implicitly based on an intuitionistically objectionable logic that assumes Aristotle's particularisation\footnote{Definition \ref{sec:A.def.1}.} is valid over $N$.

\vspace{+1ex}
\noindent However, we now show that if we assume the standard interpretation of PA is sound\footnote{In the sense of Definitions \ref{sec:A.def.9} and \ref{sec:A.def.10}.}, then PA is consistent if, and only if, it is $\omega$-consistent.

\subsubsection{Hilbert's $\omega$-Rule}
\label{sec:4.3}

To place the issue in the perspective of this paper, we consider the question:

\begin{quote}
Assuming that PA has a sound\footnote{In the sense of Definitions \ref{sec:A.def.9} and \ref{sec:A.def.10}.} interpretation over $N$, is it true that:

\textbf{Algorithmic $\omega$-Rule}: If it is proved that the PA formula $[F(x)]$ interprets as an arithmetical relation $F^{*}(x)$ that is algorithmically computable as true for any given natural number $n$, then the PA formula $[(\forall x)F(x)]$ can be admitted as an initial formula (\textit{axiom}) in PA?
\end{quote}

\noindent The significance of this query is that, as part of his program for giving mathematical reasoning a finitary foundation, Hilbert\footnote{cf.\ \cite{Hi30}, pp.485-494.} proposed an $\omega$-Rule as a finitary means of extending a Peano Arithmetic to a possible completion (i.e. to logically showing that, given any arithmetical proposition, either the proposition, or its negation, is formally provable from the axioms and rules of inference of the extended Arithmetic).

\begin{quote}
\textbf{Hilbert's $\omega$-Rule}: If it is proved that the PA formula $[F(x)]$ interprets as an arithmetical relation $F^{*}(x)$ that is true for any given natural number $n$, then the PA formula $[(\forall x)F(x)]$ can be admitted as an initial formula (\textit{axiom}) in PA.
\end{quote}

\noindent Now, in his 1931 paper---which can, not unreasonably, be seen as the outcome of a presumed attempt to validate Hilbert's $\omega$-rule---G\"{o}del introduced the concept of $\omega$-consistency\footnote{\cite{Go31}, p.23.}, from which it follows that:

\begin{lemma}
If we meta-assume Hilbert's $\omega$-rule for PA, then a consistent PA is necessarily $\omega$-consistent\footnote{However, we cannot similarly conclude from the the Algorithmic $\omega$-Rule that a consistent PA is necessarily $\omega$-consistent.}. \hfill $\Box$
\end{lemma}

\vspace{+1ex}
\noindent \textbf{Proof} If the PA formula $[F(x)]$ interprets as an arithmetical relation $F^{*}(x)$ that is true for any given natural number $n$, and the PA formula $[(\forall x)F(x)]$ can be admitted as an initial formula (\textit{axiom}) in PA, $\neg [(\forall x)F(x)]$ cannot be PA-provable if PA is consistent. The lemma follows.\hfill $\Box$

\vspace{+1ex}
\noindent Moreover, it follows from G\"{o}del's 1931 paper that one consequence of assuming Hilbert's $\omega$-Rule is that there must, then, be an undecidable arithmetical proposition\footnote{G\"{o}del constructed an arithmetical proposition $[R(x)]$ and showed that, if a Peano Arithmetic is $\omega$-consistent, then both $[(\forall x)R(x)]$ and $[\neg (\forall x)R(x)]$ are unprovable in the Arithmetic (\cite{Go31}, p.25(1), p.26(2)).}; a further consequence of which is that PA is essentially incomplete.

\vspace{+1ex}
\noindent However, since G\"{o}del's argument in this paper---from which he concludes the existence of an undecidable arithmetical proposition---is based on the weaker (i.e., weaker than assuming Hilbert's $\omega$-rule) premise that a consistent PA can be $\omega$-consistent, the question arises whether an even weaker Algorithmic $\omega$-Rule (which, prima facie, does not imply that a consistent PA is necessarily $\omega$-consistent) can yield a finitary completion for PA as sought by Hilbert, albeit for an $\omega$-\textit{in}consistent PA.

\subsubsection{Aristotle's particularisation and $\omega$-consistency}
\label{sec:4.1}

We shall now argue that these issues are related, and that placing them in an appropriate perspective requires questioning not only the persisting belief that Aristotle's 2000-year old logic of predicates---a critical component of which is Aristotle's particularisation---remains valid even when applied over an infinite domain such as $N$, but also the basis of Brouwer's denial of the Law of the Excluded Middle following his challenge of the belief in 1908\footnote{\cite{Br08}.}.

\vspace{+1ex}
\noindent Now, we have that:

\begin{lemma}
\label{sec:4.1.lem.1}
If PA is consistent but not $\omega$-consistent, then there is some PA formula $[F(x)]$ such that, under any sound\footnote{In the sense of Definitions \ref{sec:A.def.9} and \ref{sec:A.def.10}.} interpretation---say $\mathcal{I}_{PA(N,\ Sound)}$---of PA over $N$:

\begin{quote}
(i) for any given numeral $[n]$, the PA formula $[F(n)]$ interprets as true under $\mathcal{I}_{PA(N,\ Sound)}$;

(ii) the PA formula $[\neg(\forall x)F(x)]$ interprets as true under $\mathcal{I}_{PA(N,\ Sound)}$.
\end{quote}
\end{lemma}

\vspace{+1ex}
\noindent \textbf{Proof} The lemma follows from the definition of $\omega$-consistency and from Tarski's standard definitions\footnote{\cite{Ta33}; see also \cite{Ho01} for an explanatory exposition. However, for standardisation and convenience of expression, We follow the formal exposition of Tarski's definitions given in \cite{Me64}, p.50.} of the satisfaction, and truth, of the formulas of a formal system such as PA under an interpretation as detailed in Section \ref{sec:5.4.0}.\hfill $\Box$

\vspace{+1ex}
\noindent Further:

\begin{lemma}
\label{sec:4.1.lem.3}
If the interpretation $\mathcal{I}_{PA (N)}$ admits Aristotle's particularisation over $N$\footnote{As, for instance, in \cite{Me64}, pp.51-52 V(ii).}, and the PA formula $[\neg(\forall x)F(x)]$ interprets as true under $\mathcal{I}_{PA (N)}$, then there is some \textit{unspecified} PA numeral $[m]$ such that the PA formula $[F(m)]$ interprets as false under $\mathcal{I}_{PA (N)}$.
\end{lemma}

\vspace{+1ex}
\noindent \textbf{Proof} The lemma follows from Aristotle's particularisation and Tarski's standard definitions of the satisfaction, and truth, of the formulas of a formal system such as PA under an interpretation.\hfill $\Box$

\vspace{+1ex}
\noindent Hence:

\begin{lemma}
\label{sec:4.1.lem.4}
If PA is consistent and Aristotle's particularisation holds over $N$, there can be no PA formula $[F(x)]$ such that, under any sound\footnote{In the sense of Definitions \ref{sec:A.def.9} and \ref{sec:A.def.10}.} interpretation $\mathcal{I}_{PA (N,\ Sound)}$ of PA over $N$:

\begin{quote}
(i) for any given numeral $[n]$, the PA formula $[F(n)]$ interprets as true under $\mathcal{I}_{PA (N,\ Sound)}$;

(ii) the PA formula $[\neg(\forall x)F(x)]$ interprets as true under $\mathcal{I}_{PA (N,\ Sound)}$.
\end{quote}
\end{lemma}

\vspace{+1ex}
\noindent \textbf{Proof} The lemma follows from the previous two lemma.\hfill $\Box$

\vspace{+1ex}
\noindent In other words\footnote{The above argument is made explicit in view of Martin Davis' remark in \cite{Da82}, p.129, that such a proof of $\omega$-consistency may be ``\ldots open to the objection of \textit{circularity}".}:

\begin{corollary}
\label{sec:4.1.lem.5}
If PA is consistent and Aristotle's particularisation holds over $N$, then PA is $\omega$-consistent.\hfill $\Box$
\end{corollary}

\noindent It follows that:

\begin{lemma}
\label{sec:4.1.lem.6}
If Aristotle's particularisation holds over $N$, then PA is consistent if, and only if, it is $\omega$-consistent.
\end{lemma}

\vspace{+1ex}
\noindent \textbf{Proof} If PA is $\omega$-consistent then, since $[n=n]$ is PA-provable for any given PA numeral $[n]$, we cannot have that $[\neg(\forall x)(x=x)]$ is PA-provable. Since an inconsistent PA proves $[\neg(\forall x)(x=x)]$, an $\omega$-consistent PA cannot be inconsistent.\hfill $\Box$

\begin{quote}
\footnotesize{The arguments of this section thus suggest that J.\ Barkley Rosser's `extension' of G\"{o}del's argument\footnote{\cite{Ro36}.} succeeds in avoiding an explicit assumption of $\omega$-consistency \textit{only} by implicitly appealing to Aristotle's particularisation.}
\end{quote}

\noindent It further follows that:

\begin{corollary}
\label{sec:4.2.cor.1}
If PA is consistent but not $\omega$-consistent, then Aristotle's particularisation does not hold over $N$. \hfill $\Box$
\end{corollary}

\noindent As the classical, `standard', interpretation of PA---say $\mathcal{I}_{PA(N,\ Standard)}$---appeals to Aristotle's particularisation\footnote{See, for instance, \cite{Me64}, p.107 and p.52(V)(ii).}, it follows that:

\begin{corollary}
\label{sec:4.2.cor.2}
If PA is consistent but not $\omega$-consistent, then the standard interpretation $\mathcal{I}_{PA(N,\ Standard)}$ of PA is not sound\footnote{In the sense of Definitions \ref{sec:A.def.9} and \ref{sec:A.def.10}.}, and does not yield a model of PA. \hfill $\Box$
\end{corollary}

\section{Appendix B: The standard interpretation \\ $\mathcal{I}_{PA(N,\ Standard)}$ of PA over $N$}
\label{sec:7}

We define the `satisfiability' and `truth' of the formulas of PA under the standard interpretation $\mathcal{I}_{PA(N,\ Standard)}$ of PA over $N$ formally as follows:

\begin{definition}
\label{sec:7.def.1}
If $[A]$ is an atomic formula $[A(x_{1}, x_{2}, \ldots, x_{n})]$ of PA, then the denumerable sequence $(a_{1}, a_{2}, \ldots)$ in the domain $N$ of the interpretation $\mathcal{I}_{PA(N,}$ $_{Standard)}$ of PA satisfies $[A]$ if, and only if:

\begin{quote}
(i) $[A(x_{1}, x_{2}, \ldots, x_{n})]$ interprets under $\mathcal{I}_{PA(N,\ Standard)}$ as a unique relation $A^{*}(x_{1}, x_{2},$ $\ldots, x_{n})$ in $N$ for any witness $\mathcal{W}_{N}$ of $N$;

\vspace{+1ex}
(ii) for any atomic formula $[A(x_{1}, x_{2}, \ldots, x_{n})]$ of PA, and any given denumerable sequence $(b_{1}, b_{2}, \ldots)$ of $N$, there is a deterministic algorithm that computes $A^{*}(b_{1}, b_{2}, \ldots, b_{n})$ and provides objective \textit{evidence} by which any witness $\mathcal{W}_{N}$ of $N$ can \textbf{define} whether the proposition $A^{*}(b_{1}, b_{2}, \ldots, b_{n})$ holds or not in $N$;

\vspace{+1ex}
(iii) $A^{*}(a_{1}, a_{2}, \ldots, a_{n})$ holds in $N$ for any $\mathcal{W}_{N}$.
\end{quote}
\end{definition}

\noindent We inductively assign truth values of `satisfaction', `truth', and `falsity' to the compound formulas of PA under the interpretation $\mathcal{I}_{PA(N,\ Standard)}$ in terms of \textit{only} the satisfiability of the atomic formulas of PA over $N$ as follows\footnote{Compare \cite{Me64}, p.51; \cite{Mu91}.}:

\begin{definition}
\label{sec:7.def.2}
A denumerable sequence $s$ of $N$ satisfies $[\neg A]$ under $\mathcal{I}_{PA(N,\ Stan-}$ $_{dard)}$ if, and only if, $s$ does not satisfy $[A]$;
\end{definition}

\begin{definition}
\label{sec:7.def.3}
A denumerable sequence $s$ of $N$ satisfies $[A \rightarrow B]$ under $\mathcal{I}_{PA(N,}$ $_{Standard)}$ if, and only if, either it is not the case that $s$ satisfies $[A]$, or $s$ satisfies $[B]$;
\end{definition}

\begin{definition}
\label{sec:7.def.4}
A denumerable sequence $s$ of $N$ satisfies $[(\forall x_{i})A]$ under $\mathcal{I}_{PA(N,}$ $_{Standard)}$ if, and only if, given any denumerable sequence $t$ of $N$ which differs from $s$ in at most the $i$'th component, $t$ satisfies $[A]$;
\end{definition}

\begin{definition}
\label{sec:7.def.5}
A well-formed formula $[A]$ of $N$ is true\footnote{Note that this definition of `truth' is best described as `instantiational' when compared to the corresponding `algorithmic' definition of `truth' (Definition \ref{sec:8.def.5}) in Section \ref{sec:8}. The significance of the distinction between `instantiational' and `algorithmic' methods is highlighted in Section \ref{sec:9}, Appendix D.} under $\mathcal{I}_{PA(N,\ Standard)}$ if, and only if, given any denumerable sequence $t$ of $N$, $t$ satisfies $[A]$;
\end{definition}

\begin{definition}
\label{sec:7.def.6}
A well-formed formula $[A]$ of $N$ is false under $\mathcal{I}_{PA(N,\ Standard)}$ if, and only if, it is not the case that $[A]$ is true under $\mathcal{I}_{PA(N,\ Standard)}$.
\end{definition}

\section{Appendix C: The algorithmic interpretation $\mathcal{I}_{PA(N,\ Algorithmic)}$ of PA over $N$}
\label{sec:8}

We define the `satisfiability' and `truth' of the formulas of PA under the algorithmic interpretation $\mathcal{I}_{PA(N,\ Algorithmic)}$ of PA over $N$ formally as follows:

\begin{definition}
\label{sec:8.def.1}
If $[A]$ is an atomic formula $[A(x_{1}, x_{2}, \ldots, x_{n})]$ of PA, then the denumerable sequence $(a_{1}, a_{2}, \ldots)$ in the domain $N$ of the interpretation $\mathcal{I}_{PA(N,}$ $_{Algorithmic)}$ of PA satisfies $[A]$ if, and only if:

\begin{quote}
(i) $[A(x_{1}, x_{2}, \ldots, x_{n})]$ interprets under $\mathcal{I}_{PA(N,\ Algorithmic)}$ as a unique relation $A^{*}(x_{1}, x_{2},$ $\ldots, x_{n})$ in $N$ for any witness $\mathcal{W}_{N}$ of $N$;

\vspace{+1ex}
(ii) for any atomic formula $[A(x_{1}, x_{2}, \ldots, x_{n})]$ of PA, there is a deterministic algorithm that computes $[A(x_{1}, x_{2}, \ldots, x_{n})]$ and for any given denumerable sequence $(b_{1}, b_{2}, \ldots)$ of $N$, provides objective \textit{evidence} by which any witness $\mathcal{W}_{N}$ of $N$ can \textbf{define} whether the proposition $A^{*}(b_{1}, b_{2}, \ldots, b_{n})$ holds or not in $N$;

\vspace{+1ex}
(iii) $A^{*}(a_{1}, a_{2}, \ldots, a_{n})$ holds in $N$ for any $\mathcal{W}_{N}$.
\end{quote}
\end{definition}

\noindent We inductively assign truth values of `satisfaction', `truth', and `falsity' to the compound formulas of PA under the interpretation $\mathcal{I}_{PA(N,\ Algorithmic)}$ in terms of \textit{only} the satisfiability of the atomic formulas of PA over $N$ as follows\footnote{Compare \cite{Me64}, p.51; \cite{Mu91}.}:

\begin{definition}
\label{sec:8.def.2}
A denumerable sequence $s$ of $N$ satisfies $[\neg A]$ under $\mathcal{I}_{PA(N,\ Algo-}$ $_{rithmic)}$ if, and only if, $s$ does not satisfy $[A]$;
\end{definition}

\begin{definition}
\label{sec:8.def.3}
A denumerable sequence $s$ of $N$ satisfies $[A \rightarrow B]$ under $\mathcal{I}_{PA(N,}$ $_{Algorithmic)}$ if, and only if, either it is not the case that $s$ satisfies $[A]$, or $s$ satisfies $[B]$;
\end{definition}

\begin{definition}
\label{sec:8.def.4}
A denumerable sequence $s$ of $N$ satisfies $[(\forall x_{i})A]$ under $\mathcal{I}_{PA(N,}$ $_{Algorithmic)}$ if, and only if, given any denumerable sequence $t$ of $N$ which differs from $s$ in at most the $i$'th component, $t$ satisfies $[A]$.
\end{definition}

\begin{definition}
\label{sec:8.def.5}
A well-formed formula $[A]$ of $N$ is true\footnote{Algorithmically. The significance of the distinction between `instantiaional' and `algorithmic' methods is highlighted in Section \ref{sec:9}, Appendix D.} under $\mathcal{I}_{PA(N,\ Algorithmic)}$ if, and only if, given any denumerable sequence $t$ of $N$, $t$ satisfies $[A]$;
\end{definition}

\begin{definition}
\label{sec:8.def.6}
A well-formed formula $[A]$ of $N$ is false under $\mathcal{I}_{PA(N,\ Algorithmic)}$ if, and only if, it is not the case that $[A]$ is true under $\mathcal{I}_{PA(N,\ Algorithmic)}$.
\end{definition}

\section{Appendix D: The need for explicitly distinguishing between `instantiational' and `uniform' methods}
\label{sec:9}

It is significant that both Kurt G\"{o}del (initially) and Alonzo Church (subseque- ntly---possibly under the influence of G\"{o}del's disquietitude) enunciated Church's formulation of `effective computability' as a Thesis because G\"{o}del was instinctively uncomfortable with accepting it as a definition that \textit{minimally} captures the essence of `\textit{intuitive} effective computability'\footnote{See \cite{Si97}.}.

\vspace{+1ex}
\noindent G\"{o}del's reservations seem vindicated if we accept that a number-theoretic function can be effectively computable instantiationally (in the sense of being algorithmically \textit{verifiable} as envisaged in Definition \ref{sec:1.02.def.1} above), but not by a uniform method (in the sense of being algorithmically \textit{computable} as envisaged in Definition \ref{sec:1.02.def.2}).

\vspace{+1ex}
\noindent The significance of the fact (considered above in Section \ref{sec:5.4.0}) that `truth' too can be effectively decidable \textit{both} instantiationally \textit{and} by a uniform (algorithmic) method under the standard interpretation of PA is reflected in G\"{o}del's famous 1951 Gibbs lecture\footnote{\cite{Go51}.}, where he remarks: 

\begin{quote}
``I wish to point out that one may conjecture the truth of a universal proposition (for example, that I shall be able to verify a certain property for any integer given to me) and at the same time conjecture that no general proof for this fact exists. It is easy to imagine situations in which both these conjectures would be very well founded. For the first half of it, this would, for example, be the case if the proposition in question were some equation $F(n) = G(n)$ of two number-theoretical functions which could be verified up to very great numbers $n$."\footnote{Parikh's paper \cite{Pa71} can also be viewed as an attempt to investigate the consequences of expressing the essence of G\"{o}del's remarks formally.} 
\end{quote}

\noindent Such a possibility is also implicit in Turing's remarks\footnote{\cite{Tu36}, \S9(II), p.139.}:

\begin{quote}
``The computable numbers do not include all (in the ordinary sense) definable numbers. Let P be a sequence whose \textit{n}-th figure is 1 or 0 according as \textit{n} is or is not satisfactory. It is an immediate consequence of the theorem of \S8 that P is not computable. It is (so far as we know at present) possible that any assigned number of figures of P can be calculated, but not by a uniform process. When sufficiently many figures of P have been calculated, an essentially new method is necessary in order to obtain more figures." 
\end{quote}

\noindent The need for placing such a distinction on a formal basis has also been expressed explicitly on occasion\footnote{Parikh's distinction between `decidability' and `feasibility' in \cite{Pa71} also appears to echo the need for such a distinction.}. Thus, Boolos, Burgess and Jeffrey\footnote{\cite{BBJ03}, p. 37.} define a diagonal function, $d$, any value of which can be decided effectively, although there is no deterministic algorithm that can effectively compute $d$. 

\vspace{+1ex}
\noindent Now, the straightforward way of expressing this phenomenon should be to say that there are well-defined number-theoretic functions that are effectively computable instantiationally but not algorithmically. Yet, following Church and Turing, such functions are labeled as uncomputable\footnote{The issue here seems to be that, when using language to express the abstract objects of our individual, and common, mental `concept spaces', we use the word `exists' loosely in three senses, without making explicit distinctions between them (see \cite{An07c}).}!

\begin{quote}
\noindent ``According to Turing's Thesis, since $d$ is not Turing-computable, $d$ cannot be effectively computable. Why not? After all, although no Turing machine computes the function $d$, we were able to compute at least its first few values, For since, as we have noted, $f_{1} = f_{1} = f_{1} =$ the empty function we have $d(1) = d(2) = d(3) = 1$. And it may seem that we can actually compute $d(n)$ for any positive integer $n$---if we don't run out of time."\footnote{\cite{BBJ03}, p.37.}
\end{quote}

\noindent The reluctance to treat a function such as $d(n)$---or the function $\Omega(n)$ that computes the $n^{th}$ digit in the decimal expression of a Chaitin constant $\Omega$\footnote{Chaitin's Halting Probability is given by $0 < \Omega = \sum2^{-|p|} < 1$, where the summation is over all self-delimiting programs $p$ that halt, and $|p|$ is the size in bits of the halting program $p$; see \cite{Ct75}.}---as computable, on the grounds that the `time' needed to compute it increases monotonically with $n$, is curious\footnote{The incongruity of this is addressed by Parikh in \cite{Pa71}.}; the same applies to any total Turing-computable function $f(n)$\footnote{The only difference being that, in the latter case, we know there is a common `program' of constant length that will compute $f(n)$ for any given natural number $n$; in the former, we know we may need distinctly different programs for computing $f(n)$ for different values of $n$, where the length of the program will, sometime, reference $n$.}!

\section{Appendix E: No model of PA can admit a transfinite ordinal}
\label{sec:10}

Let [\(G(x)\)] denote the PA-formula:

\begin{quote}
\([x=0 \vee \neg(\forall y)\neg(x=y^{\prime})]\)
\end{quote}

\noindent Under the standard interpretation of FOL this translates, under every unrelativised interpretation of PA, as:

\begin{quote}
If \(x\) denotes an element in the domain of an unrelativised interpretation of PA, either \(x\) is 0, or \(x\) is a `successor'.
\end{quote}

\noindent Further, in every such interpretation of PA, if \(G(x)\) denotes the interpretation of [\(G(x)\)]:

\begin{quote}
(a)	\(G(0)\) is true;

(b)	If \(G(x)\) is true, then \(G(x^{\prime})\) is true.
\end{quote}

\noindent Hence, by G\"{o}del's completeness theorem:

\begin{quote}
(c)	PA proves \([G(0)]\);

(d)	PA proves \([G(x) \rightarrow G(x^{\prime})]\).

\begin{quote}
\footnotesize
\textit{G\"{o}del's Completeness Theorem}: In any first-order predicate calculus, the theorems are precisely the logically valid well-formed formulas (\textit{i.\ e.\ those that are true in every model of the calculus}).
\end{quote}
\end{quote}

\noindent Further, by Generalisation:

\begin{quote}
(e)	PA proves \([(\forall x)(G(x) \rightarrow G(x^{\prime}))]\);
\end{quote}

\noindent Hence, by Induction:

\begin{quote}
(f)	\([(\forall x)G(x)]\) is provable in PA.
\end{quote}

\noindent In other words, except 0, every element in the domain of any unrelativised interpretation of PA is a `successor'. Further, \(x\) can only be a `successor' of a unique element in any such interpretation of PA. 

\subsection{PA and Ordinal Arithmetic have no common model}
\label{sec:10.1}

Now, since Cantor's first limit ordinal, \(\omega\), is not the `successor' of any ordinal in the sense required by the PA axioms, and since there are no infinitely descending sequences of ordinals\footnote{cf.\ \cite{Me64}, p261.} in a model---if any---of set-theory, PA and Ordinal Arithmetic\footnote{cf.\ \cite{Me64}, p.187.} cannot have a common model, and so we cannot consistently extend PA to OA simply by the addition of more axioms.

\subsection{\textit{Why} PA has no set-theoretical model}
\label{sec:10.1}

We can define the usual order relation `\(<\)' in PA so that every instance of the Induction Axiom schema, such as, say:

\begin{quote}
(i) [\(F(0) \rightarrow ((\forall x)(F(x) \rightarrow F(x')) \rightarrow (\forall x)F(x))\)]
\end{quote}

\noindent yields the PA theorem:

\begin{quote}
(ii) [\(F(0) \rightarrow ((\forall x) ((\forall y)(y < x \rightarrow F(y)) \rightarrow F(x)) \rightarrow (\forall x)F(x))\)]
\end{quote}

\noindent Now, if we interpret PA without relativisation in ZF in the sense indicated by Feferman \cite{Fe92} --- i.e., numerals as finite ordinals, [\(x'\)] as [\(x \cup \left \{ x \right \}\)], etc. --- then (ii) always translates in ZF as a theorem:

\begin{quote}
(iii) [\(F(0) \rightarrow ((\forall x)((\forall y)(y \in x \rightarrow F(y)) \rightarrow F(x)) \rightarrow (\forall x)F(x))\)]
\end{quote}

\noindent However, (i) does not always translate similarly as a ZF-theorem (\textit{which is why PA and ZF can have no common model}), since the following is not necessarily provable in ZF:

\begin{quote}
(iv) [\(F(0) \rightarrow ((\forall x)(F(x) \rightarrow F(x \cup \left \{x\right \})) \rightarrow (\forall x)F(x))\)]
\end{quote}

\textit{Example}: Define [\(F(x)\)] as `[\(x \in \omega\)]'. 

\vspace{+1ex}
\noindent A significant point which emerges from the above is that we cannot appeal unrestrictedly to reasoning over transfinite ordinals when studying the foundational framework of PA.

\vspace{+1ex}
Reason: The language of PA has no constant that interprets in any model of PA as the set \textit{\textit{N}} of all natural numbers.

\vspace{+1ex}
\noindent Moreover, the preceding sections show that the Induction Axiom Schema of PA does not allow us to bypass this constraint by introducing an ``actual" (\textit{or ``completed"}) infinity disguised as an arbitrary constant - usually denoted by \(c\) or \(\infty\) - into either the language, or a putative model, of PA.

\section{Appendix F: Why the usual arguments \textit{for} a non-standard model of PA are unconvincing}
\label{sec:11}
Although we \textit{can} define a model of Arithmetic with an infinite descending sequence of elements\footnote{eg.\ \cite{BBJ03}, Section 25.\ 1, p303.}, any such model is isomorphic to the ``\textit{true arithmetic}\footnote{\cite{BBJ03}.\ p150.\ Ex.\ 12.\ 9.}" of the integers (\textit{negative plus positive}), and \textit{not} to any model of PA\footnote{\cite{BBJ03}.\ Corollary 25.\ 3, p306.}.

\vspace{+1ex}
\noindent Moreover---as we show in the next section---we cannot assume that we can consistently add a constant \(c\) to PA, along with the denumerable axioms [\(\neg (c = 0)\)], [\(\neg (c = 1)\)], [\(\neg (c = 2)\)], \ldots, since this would presume that which is sought to be proven, viz., that PA has a non-standard model.

\vspace{+1ex}
\noindent We \textit{cannot} therefore---as suggested in standard texts\footnote{eg.\ \cite{BBJ03}.\ p306; \cite{Me64}, p112, Ex.\ 2.}---apply the Compactness Theorem and the (\textit{upward}) L\"{o}wen-heim-Skolem Theorem to conclude that PA has a non-standard model.

\begin{quote}
\footnotesize
\textit{Compactness Theorem}: If every finite subset of a set of sentences has a model, then the whole set has a model\footnote{\cite{BBJ03}.\ p147.}.

\addvspace{+1ex}
\textit{Upward L\"{o}wenheim-Skolem Theorem}: Any set of sentences that has an infinite model has a non-denumerable model\footnote{\cite{BBJ03}.\ p163.}.
\end{quote}

\subsection{A formal argument \textit{for} a non-standard model of PA}

The following argument\footnote{\cite{Ln08}.} attempts to validate the above line of reasoning suggested by standard texts \textit{for} the existence of non-standard models of PA:

\begin{tabbing}
\indent 1. \= Let \(<\)\(N\) (\textit{the set of natural numbers}); \(=\) (\textit{equality}); \('\) (\textit{the successor fun-} \\ \> \textit{ction}); \(+\) (\textit{the addition function}); \( \ast \) (\textit{the product function}); \(0\) (\textit{the null} \\ \> \textit{element})\(>\) be the structure that serves to define a \textit{sound} interpretat- \\ \> ion of PA, say [\(N\)]. \\ \\

\indent 2. \> Let T[\(N\)] be the set of PA-formulas that are satisfied or true in [\(N\)]. \\ \\

\indent 3. \> The PA-provable formulas form a subset of T[\(N\)]. \\ \\

\indent 4. \> Let \( \Gamma \) be the countable set of all PA-formulas of the form \([c_{n} = (c_{n+1})']\), \\ \> where the index \(n\) is a natural number. \\ \\

\indent 5. \> Let T be the union of \( \Gamma \) and T[\(N\)]. \\ \\

\indent 6. \> T[\(N\)] plus any finite set of members of \( \Gamma \) has a model, e.g., [\(N\)] itself, \\ \> since [\(N\)] is a model of any finite descending chain of successors. \\ \\

\indent 7. \> Consequently, by Compactness, T has a model; call it \(M\). \\ \\

\indent 8. \> \(M\) has an infinite descending sequence with respect to \('\) because it is a \\ \> model of \( \Gamma \). \\ \\

\indent 9. \> Since PA is a subset of T, \(M\) is a non-standard model of PA.
\end{tabbing}

\noindent Now, if---as claimed above---[\(N\)] is a model of T[\(N\)] plus any finite set of members of \( \Gamma \), then all PA-formulas of the form \([c_{n} = (c_{n+1})']\) are PA-provable, \( \Gamma \) is a proper sub-set of the PA-provable formulas, and T is identically  T[\(N\)].

\begin{quote}
\footnotesize
The argument cannot be that some PA-formula of the form \([c_{n} = (c_{n+1})']\) is true in [\(N\)], but not PA-provable, as this would imply that PA+\([\neg (c_{n} = (c_{n+1})')]\) has a model other than [\(N\)]; in other words, it would presume that PA has a non-standard model.\footnote{The same objection applies to the usual argument found in standard texts (eg.\ [BBJ03].\ p306; \cite{Me64}, p112, Ex.\ 2) which, again, is essentially that, if PA has a non-standard model at all, then one such model is obtained by assuming we can consistently add a single non-numeral constant \(c\) to the language of PA, and the countable axioms \(c \neq 0\), \(c \neq 1\), \(c \neq 2\), \ldots to PA. However, as noted earlier, this argument too does not resolve the question of whether such assumption validly allows us to conclude that there is a non-standard model of PA in the first place.

To place this distinction in perspective, Legendre and Gauss independently conjectured in 1796 that, if \(\pi (x)\) denotes the number of primes less than \(x\), then \(\pi (x)\) is asymptotically equivalent to \(x\)/In\((x)\). Between 1848/1850, Chebyshev confirmed that if \(\pi (x)\)/\{\(x\)/In\((x)\)\} has a limit, then it must be 1. However, the crucial question of whether \(\pi (x)\)/\{\(x\)/In\((x)\)\} has a limit at all was answered in the affirmative independently by Hadamard and de la Vall\'{e}e Poussin only in 1896.}
\end{quote}

\noindent Consequently, the postulated model \(M\) of T in (7), by ``Compactness", is the model [\(N\)] that defines T[\(N\)]. However, [\(N\)] has no infinite descending sequence with respect to \('\), even though it is a model of \( \Gamma \). Hence the argument does not establish the existence of a non-standard model of PA with an infinite descending sequence with respect to the successor function \('\).

\subsection{The (\textit{upward}) Skolem-L\"{o}wenheim theorem applies only to first-order theories that admit an axiom of infinity}

We note, moreover, that the non-existence of non-standard models of PA would not contradict the (\textit{upward}) Skolem-L\"{o}wenheim theorem, since the proof of this theorem implicitly limits its applicability amongst first-order theories to those that are consistent with an axiom of infinity---in the sense that the proof implicitly requires that a constant, say \(c\), along with a denumerable set of axioms to the effect that \(c \neq 0, c \neq 1, \ldots \), can be consistently added to the theory. However, as seen in the previous section, this is not the case with PA.

\noindent \tiny{Authors postal address: 32 Agarwal House, D Road, Churchgate, Mumbai - 400 020, Maharashtra, India.\ Email: re@alixcomsi.com, anandb@vsnl.com.}

\end{document}